\documentclass[11pt]{amsart}
\usepackage{amsmath}
\usepackage{a4wide}
\usepackage[utf8]{inputenc}
\usepackage{amssymb}
\usepackage{amsopn}
\usepackage{epsfig}
\usepackage{amsfonts}
\usepackage{latexsym}
\usepackage{amsthm}
\usepackage{enumerate}
\usepackage[UKenglish]{babel}
\usepackage{verbatim}
\usepackage{color}
\usepackage{pgf}
\usepackage{tikz}

\usepackage{mathrsfs}   

\usetikzlibrary{arrows,automata}

\newtheorem{theorem}{Theorem}[section]
\newtheorem{lemma}[theorem]{Lemma}
\newtheorem{proposition}[theorem]{Proposition}
\newtheorem{corollary}[theorem]{Corollary}

\theoremstyle{definition}

\theoremstyle{remark}
\newtheorem{remark}[theorem]{Remark}

\numberwithin{equation}{section}

\newcommand{\R}{\ensuremath{\mathbb{R}}}
\newcommand{\N}{\ensuremath{\mathbb{N}}}

\renewcommand{\c}{ {\mathbf{c}}}
\renewcommand{\d}{ {\mathbf{d}}}
\newcommand{\T}{\ensuremath{\mathbf{T}}}

\newcommand{\set}[1]{\left\{#1\right\}}
\newcommand{\la}{\lambda}
\newcommand{\ga}{\gamma}
\newcommand{\Ga}{\Gamma}

\newcommand{\f}{\infty}
\newcommand{\de}{\alpha^*}

\newcommand{\al}{\alpha}
\newcommand{\lle}{\preccurlyeq}
\newcommand{\lge}{\succcurlyeq}

\newcommand{\La}{\Lambda}
\newcommand{\F}{\mathcal{F}}

\renewcommand{\de}{\delta}

\begin{document}

\title{Intersections of middle-$\alpha$  Cantor sets with a fixed translation}

\author[Y.Huang]{Yan Huang}
\address[Y. Huang]{College of Mathematics and Statistics, Chongqing University, 401331, Chongqing, P.R.China}
\email{{yanhuangyh@126.com}}

\author[D.Kong]{Derong Kong}
\address[D. Kong]{College of Mathematics and Statistics, Chongqing University, 401331, Chongqing, P.R.China}
\email{derongkong@126.com}




\begin{abstract}
For $\la\in(0,1/3]$ let $C_\la$ be the middle-$(1-2\la)$ Cantor set in $\R$. Given $t\in[-1,1]$, excluding the trivial {case}  we show that
{\[
\La(t):=\set{\la\in(0,1/3]: C_\la\cap(C_\la+t)\ne\emptyset}
\]
is a topological Cantor set with zero Lebesgue measure and full Hausdorff dimension.} In particular, we calculate the local dimension of $\La(t)$, which reveals a dimensional variation principle. Furthermore, for any $\beta\in[0,1]$ we show that the level set
\[
\La_\beta(t):=\set{\la\in\La(t): \dim_H(C_\la\cap(C_\la+t))=\dim_P(C_\la\cap(C_\la+t))=\beta\frac{\log 2}{-\log \la}}
\]
has {equal} Hausdorff and packing dimension  {$(-\beta\log\beta-(1-\beta)\log\frac{1-\beta}{2})/\log 3$}. We also show that the set of $\la\in\La(t)$ for which $\dim_H(C_\la\cap(C_\la+t))\ne\dim_P(C_\la\cap(C_\la+t))$ has full Hausdorff dimension.
\end{abstract}

\maketitle

\section{Introduction}
Intersections of Cantor sets on the real line {appear} in the setting of homoclinic bifurcations in dynamical systems (cf.~\cite{Moreira-1996}). Moreira and Yoccoz \cite{Moreira-Yoccoz-2001} studied the stable intersections of regular Cantor sets, and gave an affirmative answer to Palis' conjecture. Intersections of Cantor sets also {appear} in number theory. Hall \cite{Hall-1947} proved that any real number can be written as the sum of two numbers whose continued fractional coefficients are at most $4$, and from this it follows that the Lagrange spectrum  contains a whole half-line. {Note} that the middle-$\alpha$ Cantor set is an affine Cantor set, which minimizes the Hausdorff dimension with a given thickness (cf.~\cite{Kraft-1992, Palis_Takens_1993}). Motivated by the above works there is a great interest in the study of intersections of middle-$\alpha$ Cantor set with its translations. {Kraft \cite{Kraft-1994} gave a complete description when the intersection of middle-$\alpha$ Cantor set with its translation is a single point. Li and Xiao \cite{Li-Xiao-1998} calculated the Hausdorff and packing dimensions of the intersection of middle-$\alpha$ Cantor set with its translation  for $\alpha\ge 1/3$. When $\alpha\in(0,1/3)$,  Zou, Lu and Li \cite{Zou-Lu-Li-2008} determined when the intersection of middle-$\alpha$ Cantor set with its translation is a self-similar set under the condition that the translation has a unique coding. Recently, Baker and the second author \cite{Baker-Kong-2017} proved that for $\al\in(0,1/3)$ it is possible {that} the intersection of middle-$\alpha$ Cantor set with its translation contains only Liouville numbers.}

{For} $\la\in(0,{1/2})$ let $C_\la$ be the middle-$\alpha$ Cantor set with $\alpha=1-2\la$. Then $C_\la$ is {a  \emph{self-similar set}} generated by the {\emph{iterated function system} (IFS):
$
\set{g_{i}(x)=\la x+i(1-\la): i=0, 1}.
$}
In other words,  $C_\la$ is the unique nonempty compact set satisfying $C_\la=g_0(C_\la)\cup g_1(C_\la)$. {So,
\begin{equation}\label{eq:C-lambda}
C_\la=\set{(1-\la)\sum_{n=1}^\f i_n\la^{n-1}: i_n\in\{0,1\}~\forall n\ge 1}.
\end{equation}
Observe that}
$C_\la\cap(C_\la+t)\ne\emptyset$ if and only if $t\in C_\la-C_\la$.
By (\ref{eq:C-lambda}) it follows that
\begin{equation}\label{eq:Translation}
 E_\la:=C_\la-C_\la=\set{(1-\la)\sum_{n=1}^\infty {i_n}\la^{n-1}:{i_n}\in\{-1,0,1\}~\forall n\ge 1}.
\end{equation}
  Clearly, if $\la\in(0,1/3)$, then $E_\la$ is a Cantor set having zero Lebesgue measure. {And}  each $t\in E_\la$ has a unique \emph{coding}  ${(i_n)} \in\set{-1,0,1}^\N$ such that
$
 t=(1-\la)\sum_{n=1}^\f {i_n}\la^{n-1}.
 $
  Li and Xiao \cite[{Theorem 3.4}]{Li-Xiao-1998} gave the Hausdorff and packing dimensions of $C_\la\cap(C_\la+t)$:
  {
 \begin{equation}\label{eq:dim-intersection}
 \begin{split}
 &\dim_H(C_\la\cap(C_\la+t))=\frac{\log 2}{-\log\la}\,\liminf_{n\to\infty}\frac{\#\{1\leq j \leq n: i_j=0\}}{n},\\
 &\dim_P(C_\la\cap(C_\la+t))=\frac{\log 2}{-\log\la}\,\limsup_{n\to\infty}\frac{\#\{1\leq j \leq n: i_j=0\}}{n},
 \end{split}
 \end{equation}
   where $\#A$ denotes the cardinality of a set $A$.}
   If $\la=1/3$, then   except for  a countable set of   {points} in $E_\la$  having two different codings, all other   {$t$s in $E_\la$} have a unique coding; and the dimension formulae (\ref{eq:dim-intersection}) still {hold.} If $\la\in(1/3,1/2)$, then $E_\la=[-1,1]$, and it is well known that Lebesgue almost every {$t \in [-1,1]$} has a continuum of codings (cf.~\cite{Sidorov_2003}). In this case the dimension formulae   (\ref{eq:dim-intersection}) {fail} for typical $t\in[-1,1]$, but we still have the dimension formulae (\ref{eq:dim-intersection}) if $t$ has a unique coding (see \cite{Kong_Li_Dekking_2010}).

\begin{figure}[h!]
  \centering
  \includegraphics[width=8cm]{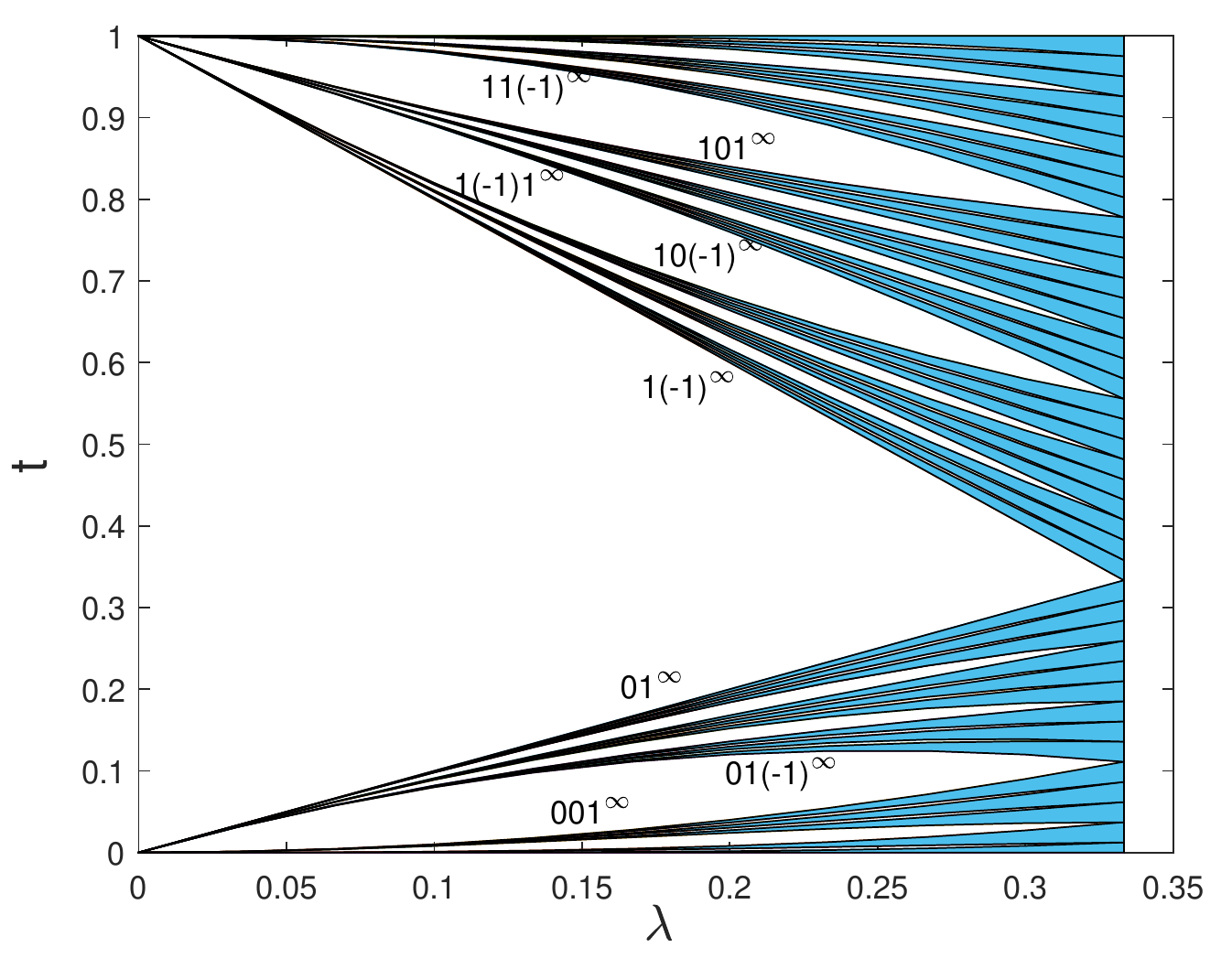}\\
  \caption{The {fourth} level approximation of the master set $\Gamma$ which consists of all vectors $(\la, t)\in(0,1/3]\times[0,1]$ such that $C_\la\cap(C_\la+t)\ne\emptyset$. Each curve corresponds to a unique coding $(i_n)\in\set{-1,0,1}^\N$ via the equation $t=(1-\la)\sum_{n=1}^\f i_n\la^{n-1}$. {For any $(\lambda,t)\in\Gamma$ the} vertical fiber is $E_\lambda\cap [0,1]$, and the horizontal fiber is $\Lambda(t)$. }\label{Fig:0-1}
\end{figure}

 Let $\tilde\Ga:=\set{(\la, t)\in(0,1/2)\times[-1,1]: C_\la\cap(C_\la+t)\ne\emptyset}$ be the master set defined by the intersections of Cantor sets. Since for any $\la\in(1/3,1/2)$ and $t\in[-1,1]$ the intersection $C_\la\cap(C_\la+t)$ is always non-empty, and $C_\la\cap(C_\la+t)\ne\emptyset$ if and only if $C_\la\cap(C_\la-t)\ne\emptyset$, it is then interesting to consider the subset (see Figure \ref{Fig:0-1})
 \[
 \Ga:=\set{(\la, t)\in(0,\frac{1}{3}]\times [0,1]: C_\la\cap(C_\la+t)\ne\emptyset}=\set{(\la, t){\in(0,\frac{1}{3}]\times[0,1]}: t\in E_\la},
 \]
where the second equality {holds} because $C_\la\cap(C_\la+t)\ne\emptyset$ if and only $t\in C_\la-C_\la=E_\la$.
{Observe} that for $\la\in(0,1/3]$ the vertical fiber $\Ga(\la)=\set{t\in[0,1]: t\in E_\la}=E_{\la}\cap[0,1]$ is a Cantor set having Hausdorff dimension $-\log 3/\log\la$. Since $E_\la$ is a self-similar set, a lot {is known} about this vertical fiber $\Ga(\la)$ (cf.~\cite{Hutchinson_1981}). On the other hand, {for any $t\in[0,1]$}  let
\begin{equation}\label{eq:def-Lambda-t}
\La(t):=\set{\la\in(0,\frac{1}{3}]:C_\la\cap(C_\la+t)\ne\emptyset}=\set{\la\in(0,\frac{1}{3}]:t\in E_\la}
\end{equation}
be a horizontal fiber of $\Ga$.
By (\ref{eq:Translation}) it follows that  {$\La(0)=\La(1)=(0,1/3]$}, and $\La(1/3)=\{1/3\}$. So {it is interesting to study $\La(t)$ for  $t\in(0,1)\setminus\{1/3\}$.

 Our first result shows that  $\La(t)$ is   a \emph{topological Cantor set}, which is a non-empty compact set having neither interior nor isolated points.}

\begin{theorem}\label{thm: La(t)}
Let $t\in(0,1)\setminus \{1/3\}$.
\begin{enumerate}
  [{\rm(i)}]
  \item The set $\La(t)$ is a topological Cantor set with ${\min}\La(t)=\min\{t,\frac{1-t}{2}\}$ and ${\max}\La(t)=1/3$;

  \item $\La(t)$ is a Lebesgue null set having full Hausdorff dimension;

  \item For any $\la\in\La(t)$ we have
  \[
  \lim_{\de\to 0^+}\dim_H(\La(t)\cap(\la-\de, \la+\de))=\frac{\log 3}{-\log\la}.
  \]
\end{enumerate}
\end{theorem}
\begin{remark}
  Note that $E_\la$ is a self-similar set having Hausdorff dimension $-\log 3/\log\la$. Then Theorem \ref{thm: La(t)} (iii) shows that the local dimension of $\La(t)$ at $\la$ is the same as the local dimension of $E_\la$ at $t$. In view of Figure \ref{Fig:0-1} it follows that the local dimensions of $\Gamma$ at any point $(\la, t){\in\Gamma}$ through the horizontal and the vertical fibers are the same{; this} can be viewed as a dimensional `variation principle' for the set $\Gamma$.
\end{remark}
{Note by Theorem \ref{thm: La(t)} (i)} that $\Lambda(t)$ is a topological Cantor set. Then it can be obtained by successively removing a sequence of open intervals from the closed interval {$[\min\Lambda(t),1/3]$}.
{A} geometrical construction of $\La(1/2)$ {is} plotted in Figure \ref{fig:1} (left).
By {Theorem \ref{thm: La(t)} {(iii)}}   it follows that
\[
\dim_H(\La(t)\cap(0,\la])=\sup_{\ga\in\La(t)\cap(0,\la)}\frac{\log 3}{-\log\ga}\qquad \forall \la\in(0,1/3].
\]
{Therefore}, the dimension function $\psi_t: \la\mapsto\dim_H(\La(t)\cap(0,\la])$, which describes the distribution of $\La(t)$, is a \emph{C\'adl\'ag function} (see {Figure \ref{fig:1}, right}).    It is locally constant almost everywhere, left continuous   with right-hand limits  everywhere, and thus it has countably {infinitely} many discontinuities; and it has no downward jumps.

\begin{center}
\begin{figure}[h!]
\begin{tikzpicture}[xscale=65,yscale=80,axis/.style={very thick, ->}, important line/.style={thick}, dashed line/.style={dashed, thin},
    pile/.style={thick, ->, >=stealth', shorten <=2pt, shorten>=2pt},every node/.style={color=black} ]

    \draw[important line] ({1/4}, 0)--({1/3}, 0);
     \node[] at({1/4}, 0.005){$\frac{1}{4}$};  \node[] at({1/3}, 0.005){$\frac{1}{3}$};

     \draw[important line] ({1/4}, -0.02)--({0.2696}, -0.02);  \draw[important line] ({0.2779}, -0.02)--({0.3154}, -0.02);
     \draw[important line] ({0.3194}, -0.02)--({1/3}, -0.02);
     \node[scale=0.6] at({0.2666}, -0.015){\scriptsize{$1(-1)^21^\f$}};  \node[scale=0.6] at({0.2809}, -0.025){\scriptsize{$1(-1)0(-1)^\f$}};
     \node[scale=0.6] at({0.3124}, -0.025){\scriptsize{$1(-1)01^\f$}}; \node[scale=0.6] at({0.3264}, -0.015){\scriptsize{$1(-1)1(-1)^\f$}};

     \draw[important line] ({1/4}, -0.04)--({0.2542}, -0.04);\draw[important line]({0.2563},-0.04)--({0.2612}, -0.04);
     \draw[important line]({0.2635},-0.04)--({0.2696}, -0.04);

     \draw[important line] ({0.2779}, -0.04)--(0.2859, -0.04);\draw[important line] ({0.2880}, -0.04)--({0.2985}, -0.04);
     \draw[important line]({0.3006},-0.04)--({0.3154}, -0.04);

     \draw[important line]({0.3194},-0.04)--({1/3}, -0.04);

    \draw[important line] ({1/4}, -0.06)--(0.2510, -0.06); \draw[important line] (0.2515, -0.06)--(0.2525, -0.06);
    \draw[important line](0.2531,-0.06)--(0.2542, -0.06);
    \draw[important line](0.2563,-0.06)--({0.2575}, -0.06);\draw[important line] ({0.2580}, -0.06)--(0.2593, -0.06);
    \draw[important line] (0.2599, -0.06)--(0.2612, -0.06);
    \draw[important line] ({0.2635}, -0.06)--(0.2650, -0.06);\draw[important line] (0.2656, -0.06)--({0.2672}, -0.06);
    \draw[important line] ({0.2678}, -0.06)--(0.2696, -0.06);

    \draw[important line] ({0.2779}, -0.06)--(0.2799, -0.06); \draw[important line] (0.2805, -0.06)--(0.2828, -0.06);
    \draw[important line](0.2834,-0.06)--(0.2859, -0.06);
    \draw[important line](0.2880,-0.06)--({0.2908}, -0.06);\draw[important line] ({0.2914}, -0.06)--(0.2945, -0.06);
    \draw[important line] (0.2951, -0.06)--(0.2985, -0.06);
    \draw[important line] ({0.3006}, -0.06)--(0.3045, -0.06);\draw[important line] (0.3051, -0.06)--({0.3096}, -0.06);
    \draw[important line] (0.3101, -0.06)--(0.3154, -0.06);

    \draw[important line] ({0.3194}, -0.06)--(0.3249, -0.06); \draw[important line] (0.3251, -0.06)--(0.3315, -0.06);
    \draw[important line](0.3316,-0.06)--({1/3}, -0.06);
    \draw[important line](0.3333,-0.065)--(0.3333, -0.065);

\end{tikzpicture}
\quad
\begin{tikzpicture}[xscale=65,yscale=24]
\draw[<->]  (.25,1.01) node[scale=1,anchor=south]{$\psi_{\frac{1}{2}}(\lambda)$}--(.25,.79) node[scale=0.6,anchor=east]{$0.79$} node[scale=0.6,anchor=north]{$0.25$} -- (.335,.79) node[scale=1,anchor=west]{$\lambda$};
\draw[dashed] (0.2696,0.787) node[scale=0.6,anchor=north]{$0.269$} --(0.2696,0.8381);
\draw (0.2800,0.787) node[scale=0.6,anchor=north]{$0.278$};
\draw[dashed] (0.2779,0.787)  --(0.2779,0.8580);

\draw (0.3130,0.787) node[scale=0.6,anchor=north]{$0.315$};
\draw[dashed] (0.3154,0.787)  --(0.3154,0.9521);
\draw (0.3219,0.787) node[scale=0.6,anchor=north]{$0.319$};
\draw[dashed] (0.3194,0.787)  --(0.3194,0.9626);

\draw (0.3333,0.787) node[scale=0.6,anchor=north]{$\frac{1}{3}$};
\draw[dashed] (0.3333,0.787)  --(0.3333,1);
\draw (0.25,1) node[scale=0.6,anchor=east]{$1$};
\draw[dashed] (0.25,1)  --(0.3333,1);

\draw[thick](.25,.7925)--(.2510,.7948)--(.2515,.7948)
(.2515,.7959)--(.2525,.7982)--(.2531,.7982)
(.2531,.7996)--(.2542,.8021)--(.2563,.8021)
(.2563,.8070)--(.2575,.8097)--(.2580,.8097)
(.2580,.8109)--(.2593,.8139)--(.2599,.8139)
(.2599,.8153)--(.2612,.8184)--(.2635,.8184)
(.2635,.8237)--(.2650,.8273)--(.2656,.8273)
(.2656,.8287)--(.2672,.8324)--(.2678,.8324)
(.2678,.8339)--(.2696,.8381)--(.2779,.8381)
(.2779,.8580)--(.2799,.8628)--(.2805,.8628)
(.2805,.8642)--(.2828,.8698)--(.2834,.8698)
(.2834,.8713)--(.2859,.8774)--(.2880,.8774)
(.2880,.8826)--(.2908,.8895)--(.2914,.8895)
(.2914,.8910)--(.2945,.8987)--(.2951,.8987)
(.2951,.9002)--(.2985,.9087)--(.3006,.9087)
(.3006,.9140)--(.3045,.9239)--(.3051,.9239)
(.3051,.9254)--(.3096,.9370)--(.3101,.9370)
(.3101,.9383)--(.3154,.9521)--(.3194,.9521)
(.3194,.9626)--(.3249,.9772)--(.3251,.9772)
(.3251,.9777)--(.3315,.9950)--(.3316,.9950)
(.3316,.9953)--(.3333,1);

\end{tikzpicture}
\caption{Left: {the  geometrical construction of $\La(1/2)$}. Right: the graph of $\psi_{1/2}: \lambda\mapsto\dim_H(\Lambda(1/2)\cap(0,\lambda])$ for $\lambda\in(\min\Lambda(1/2),1/3]=(1/4,1/3]$.   }
\label{fig:1}
\end{figure}
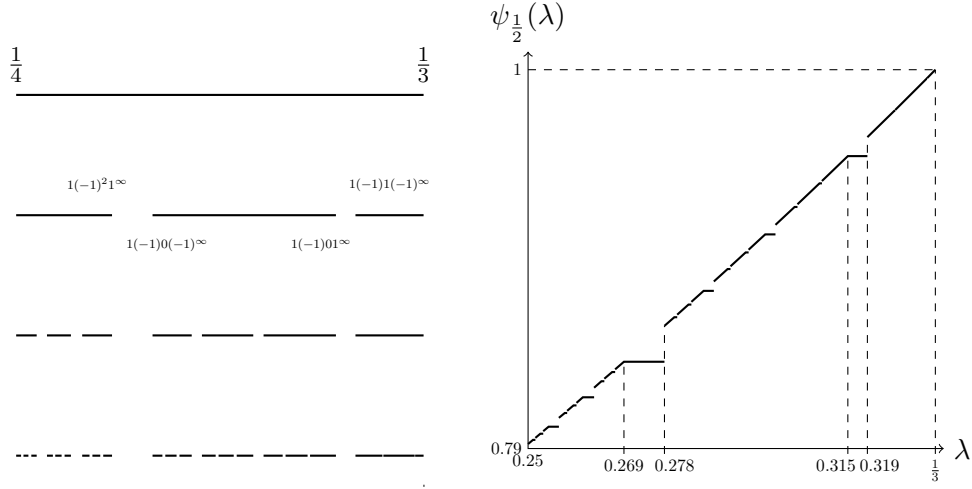
\end{center}

{Observe}  by (\ref{eq:dim-intersection}) that the dimension of $C_\la\cap(C_\la+t)$ is determined by the frequency of digit zero in the coding of $t$ {in base $\lambda$}. Then for $\la\in(0,1/3]$ the {vertical} fiber  $\Ga(\la)$ can be partitioned {as}
\[\Ga(\la)=\Ga_{not}(\la)\cup\bigcup_{\beta\in[0,1]}\Ga_\beta(\la),\]
where
\begin{align*}
\Ga_{not}(\la)&:=\set{t\in[0,1]: \dim_H (C_\la\cap(C_\la+t))\ne\dim_P (C_\la\cap(C_\la+t))},\\
\Ga_\beta(\la)&:=\set{t\in[0,1]: \dim_H (C_\la\cap(C_\la+t))=\dim_P (C_\la\cap(C_\la+t))=\beta\frac{\log 2}{-\log\la}}.
\end{align*}
Li and Xiao \cite[{Theorem 4.3}]{Li-Xiao-1998}  showed that
\[
\dim_H\Ga_\beta(\la)=\dim_P\Ga_\beta(\la)=\frac{h(\frac{1-\beta}{2}, \beta,\frac{1-\beta}{2})}{-\log\la},
\]
where {for a  probability vector $(p_1, p_2, p_3)$}
\begin{equation}\label{eq:entropy}
{h(p_1,p_2,p_3):=-\sum_{i=1}^3 p_i\log p_i}.
\end{equation}
{Here we adopt the convention $0\log0=0$.} Furthermore, an application of \cite{Barreira-Schmeling-2000} gives that
{$\dim_H\Ga_{not}(\la)=-\log 3/\log \la$}.

{Inspired by the works of \cite{Barreira-Schmeling-2000} and \cite{Li-Xiao-1998} we consider the level sets {of the horizontal fiber} $\La(t)$. Given $t\in(0,1)\setminus\set{1/3}$, for $\beta\in[0,1]$ let
\[\Lambda_\beta(t):=\set{\lambda\in\La(t): \dim_H \big(C_\lambda\cap(C_\lambda+t)\big)=\dim_P \big(C_\lambda\cap(C_\lambda+t)\big)=\beta\;\frac{\log 2}{-\log \la}    }.\]
 Then the horizontal fiber $\Lambda(t)$ can be partitioned {as}
\[
\La(t)=\La_{not}(t)\cup\bigcup_{\beta\in[0,1]}\La_\beta(t),
\]
where
\begin{align*}
\La_{not}(t)&:=\set{\la\in\La(t): \dim_H \big(C_\lambda\cap(C_\lambda+t)\big)\ne\dim_P \big(C_\lambda\cap(C_\lambda+t)\big)}.
\end{align*}

\begin{theorem}
  \label{th:La(t)-alpha}
  For any $t\in(0,1)\setminus\set{1/3}$ and $\beta\in[0,1]$  the sets $\La_{not}(t)$  and $\La_\beta(t)$ are both dense in $\La(t)$. Furthermore,
  \[
\dim_H\La_{not}(t)=1\quad\textrm{and}\quad \dim_H\La_\beta(t)=\dim_P\La_\beta(t)=\frac{h(\frac{1-\beta}{2}, \beta,\frac{1-\beta}{2})}{\log 3}.
\]
\end{theorem}
}

\begin{remark}
  Our proof of Theorem \ref{th:La(t)-alpha} can be adapted to {showing} that for any $t\in(0,1)\setminus\set{1/3}$
  \[\La_{fin}(t)={\set{\la\in\La(t): \#C_\la\cap(C_\la+t)<+\f}}\]
  has the same Hausdorff dimension as $\La_0(t)$, that is $\dim_H\La_{fin}(t)=\log 2/\log 3$.
\end{remark}

The rest of the paper is organized as follows. In the next section we define {the} coding map $\Phi_t$ which maps each $\la\in\La(t)$ to its   coding of $t$ in base $\la$, and show that $\Phi_t$ is continuous and piecewise monotonic in $\La(t)$. Based on this {map {$\Phi_t$}} we show in Section \ref{sec:topology-La(t)}   that $\La(t)$ is a {topological} Cantor set, and it has full Hausdorff dimension. In Section \ref{sec:local-dimension} we calculate the local dimension of $\La(t)$, and prove Theorem \ref{thm: La(t)}.  Motivated by the {works} of non-normal numbers we show in Section \ref{sec:level-set-not}   that $\La_{not}(t)$ is a dense subset of $\La(t)$ with full Hausdorff dimension. Finally, in Section \ref{sec:Proof-theorem2} we calculate the dimension of $\La_\beta(t)$ and prove Theorem {\ref{th:La(t)-alpha}}.

\section{Preliminaries}\label{sec:preliminaries}
In this section we will define a map {$\Phi_t$} which maps $\Lambda(t)$ to the symbolic space $\{-1,0,1\}^{\mathbb{N}}$.
First we recall some terminology from symbolic dynamics (cf.~\cite{Lind_Marcus_1995}). Let $\set{-1,0, 1}^\N$ be the set of all infinite sequences over the {\emph{alphabet}} $\set{-1,0,1}$. By a \emph{word} we mean a finite string of digits over $\set{-1,0,1}$. Denote by $\set{-1,0,1}^*$ the set of all finite words including the empty word $\epsilon$. For two words $\mathbf c=c_1\ldots c_m$ and $\mathbf d=d_1\ldots d_n$ we write $\mathbf{cd}=c_1\ldots c_md_1\ldots d_n$ for their concatenation. In particular, for any $k\in\N$ we denote by $\mathbf c^k$ the $k$-fold concatenation of $\mathbf c$ with itself, and by $\mathbf c^\f$ the periodic sequence which is obtained by the infinite concatenation of $\mathbf c$ with itself. For a word {$\mathbf c=c_1\ldots c_m\in\set{-1,0,1}^*$ with $m\ge 1$, if} $c_m<1$ we write $\mathbf c^+:=c_1\ldots c_{m-1}(c_m+1)$; and if $c_m>-1$, we write $\mathbf c^-:=c_1\ldots c_{m-1}(c_m-1)$. {Therefore, $\mathbf c^+$ and $\mathbf c^-$ are both words over the {{alphabet}} $\set{-1,0,1}$.}
Throughout the paper we will use lexicographical order $\prec, \lle, \succ$ or $\lge$ between sequences in $\set{-1,0,1}^\N$. For example, we say  $(i_n)\succ (j_n)$ if $i_1>j_1$, or there exists $n\in\N$ such that $i_1\ldots i_n=j_1\ldots j_n$ and $i_{n+1}>j_{n+1}$. And we write $(i_n)\lge (j_n)$ if $(i_n)\succ (j_n)$ or $(i_n)=(j_n)$. Similarly, we say $(i_n)\prec (j_n)$ if $(j_n)\succ (i_n)$, and say $(i_n)\lle (j_n)$ if $(j_n)\lge (i_n)$.
For two infinite sequences $\c=(c_n), \d=(d_n)\in\set{-1,0,1}^\N$ with $\c\prec \d$ we write
\[(\c, \d):=\set{(i_n)\in\set{-1,0,1}^\N: \c\prec (i_n)\prec \d},\quad [\c,\d]:=\set{(i_n)\in\set{-1,0,1}^\N: \c\lle(i_n)\lle \d}.\]
Similarly, we {set}
\[(\c,\d]:=\set{(i_n)\in\set{-1,0,1}^\N: \c\prec(i_n)\lle \d},\quad [\c,\d):=\set{(i_n)\in\set{-1,0.1}^\N:\c\lle(i_n)\prec \d}.\]

{For $\la\in(0,1/3]$ we note by (\ref{eq:Translation}) that $E_\la$ is a self-similar set generated by the IFS {$\set{g_{i}(x)=\la x+i(1-\la): i=-1, 0, 1}$}. This induces a {map} $\pi_\la$ defined by
\[
\pi_\la:\set{-1,0,1}^{\mathbb{N}}\to E_\la;\quad  (i_n)\mapsto (1-\la)\sum\limits_{n=1}^\infty i_n\la^{n-1}.
\]
{It is clear that {the map $\pi_\la$ is bijective if $\la\in(0,1/3)$; and  the map $\pi_\la$ is bijective up to a countable set if $\la=1/3$}.}
 Now, based on  $\pi_\la$ we define the master  map $\Pi$ by
\begin{equation}\label{eq:map-Pi}
\Pi :\{-1,0,1\}^{\mathbb{N}} \times (0,1/3]\to [-1,1];\quad ((i_n),\lambda)\mapsto \pi_\la((i_n))=(1-\lambda)\sum_{n=1}^\infty i_n \lambda^{n-1}.
\end{equation}
{Note that the symbolic space $\set{-1,0,1}^\N$ becomes} a compact metric space under the metric $\rho$  defined by
\[
\rho((i_n), (j_n))=3^{1-\inf\set{n\ge 1: i_n\ne j_n}}.
\]

First we show that $\Pi$ is continuous under the product topology induced by the metric $\rho$ on $\set{-1,0,1}^\N$ and the Euclidean metric $|\cdot|$ on $\R$.
\begin{lemma}
  \label{lem:continuity-Pi}
  The map $\Pi: \set{-1,0,1}^\N\times(0,1/3] \to [-1,1]$ is continuous and onto.
\end{lemma}
\begin{proof}
  Note  that $\Pi(\set{-1,0,1}^\N\times\set{1/3})=[-1,1]$. {It suffices} to prove the  continuity of $\Pi$, {which}  follows from the following observation: for any two $((i_n), \la_1), ((j_n), \la_2)\in\set{-1,0,1}^\N\times (0,1/3]$ with $k=\inf\set{n\ge 1: i_n\ne j_n}$, by (\ref{eq:map-Pi}) it follows that
  \begin{align*}
|\Pi((i_n), \la_1)-\Pi((j_n), \la_2)|&\le |\Pi((i_n), \la_1)-\Pi((i_n), \la_2)|+|\Pi((i_n), \la_2)-\Pi((j_n), \la_2)|\\
&\le \sup_{\la\in(0,1/3]}|\Pi_2((i_n),\la)|\cdot|\la_1-\la_2|+\left|(1-\la_2)\sum_{n=1}^\f(i_n-j_n)\la_2^{n-1}\right|\\
&\le 2|\la_1-\la_2|+\left|(1-\la_2)\sum_{n=k}^\f 2\la_2^{n-1}\right|\\
&\le 2|\la_1-\la_2|+2 \rho((i_n), (j_n)),
  \end{align*}
  where the second inequality follows by the mean value theorem that
   \begin{equation}
    \label{eq:derivative}
  \Pi_2((i_n),\la):=\frac{\partial \Pi((i_n), \la)}{\partial\la}=\sum_{n=1}^\f i_n n\left(\frac{n-1}{n}-\la\right)\la^{n-2},
  \end{equation} and the third inequality follows by
 $
   \left|\Pi_2((i_n),\la)\right|  \le 1+\sum_{n=2}^\f n\left(\frac{n-1}{n}-\la\right)\la^{n-2}=2.
$
\end{proof}

Next we show that  $\Pi$ is monotonic {in} its first variable.}

\begin{lemma}\label{lem:monotonicity-in}
{Given} $\la\in(0,1/3]$, the function $\Pi(\cdot,\la)$ is increasing {in $\set{-1,0,1}^\N$} with respect to the lexicographical order. In particular, if $\la\in(0,1/3)$, the function $\Pi(\cdot,\la)$ is strictly increasing.
\end{lemma}
\begin{proof} Take $(i_n),(j_n)\in\{-1,0,1\}^{\mathbb{N}}$ with $(i_n){\prec}(j_n)$. Then there exists {$m \in \mathbb{N}$}  such that $i_n = j_n$ for all $n<m$, and $i_m<j_m$.   By (\ref{eq:map-Pi}) this implies
\begin{align*}
\Pi((j_n),\la)-\Pi((i_n),\la)&=(1-\lambda)\sum_{n=m}^\infty(j_n-i_n)\lambda^{n-1}\\
&\ge (1-\lambda)\lambda^{m-1}- (1-\lambda)\sum_{n=m+1}^\infty 2\lambda^{n-1}\\
&=(1-\lambda)\lambda^{m-1}-2\lambda^m\ge0,
\end{align*}
{where the last inequality follows by   $\la\in(0, 1/3]$; and this   inequality is strict if $\la\in(0,1/3)$.}
\end{proof}
\begin{figure}[h!]
\begin{center}
\includegraphics[width=7cm]{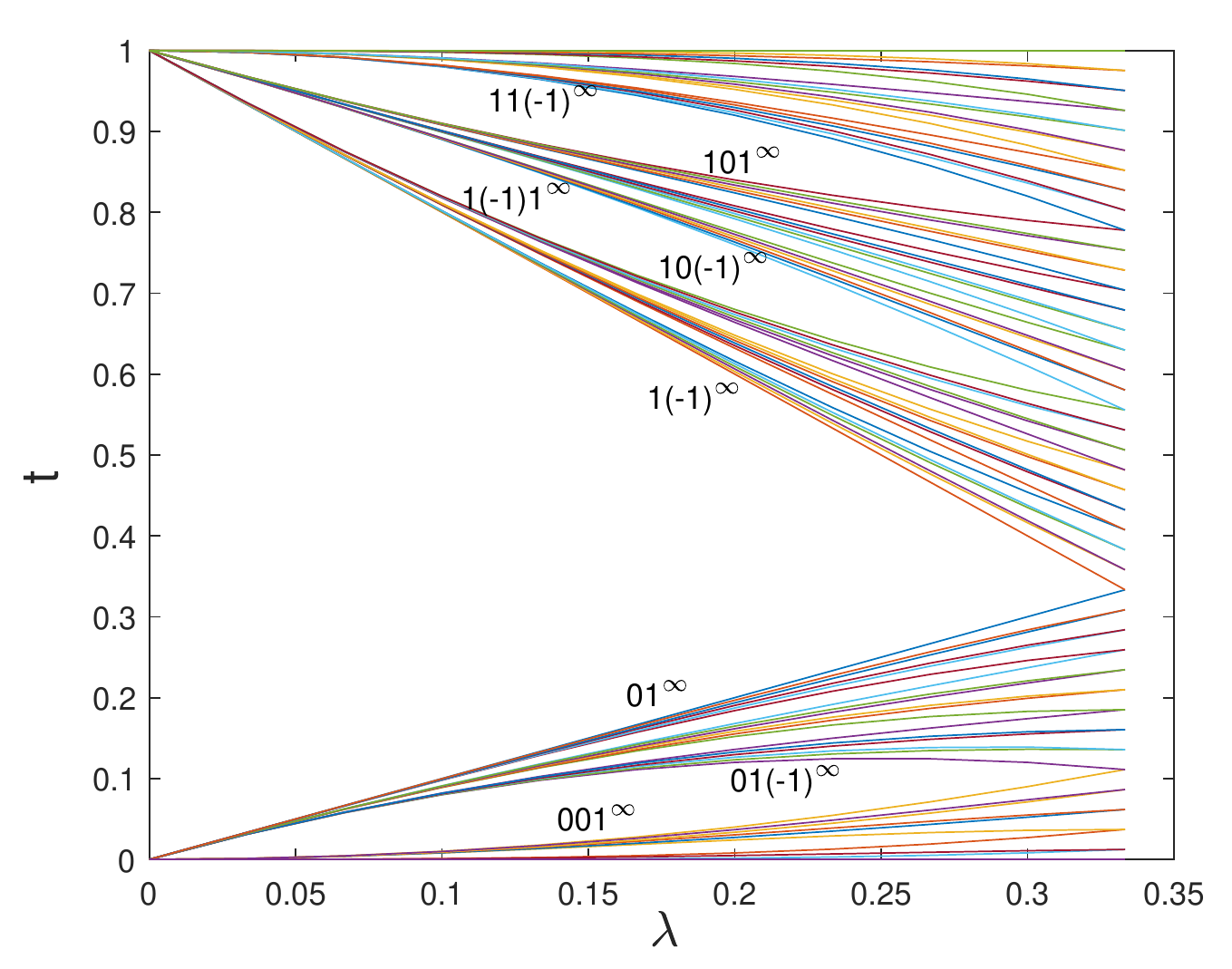}\quad
\begin{tikzpicture}[
    scale=6,
    axis/.style={very thick, ->},
    important line/.style={thick},
    dashed line/.style={dashed, thin},
    pile/.style={thick, ->, >=stealth', shorten <=2pt, shorten
    >=2pt},
    every node/.style={color=black}
    ]

    \draw[important line] (0, 0)--(1.2, 0);

    \draw[important line] (0, -0.02)--(0, 0.02);\node[scale=0.6] at(0, -0.08){$0^\f$};\node[scale=0.6] at(0, 0.08){$0$};
    \draw[-latex] (0.05, 0.05)--(0.25, 0.15);

       \draw[important line] (0.3, -0.02)--(0.3, 0.02);\node[scale=0.6] at(0.3, -0.08){$001^\f$};\node[scale=0.6] at(0.3, 0.08){$\frac{1}{9}$};
              \draw[-latex] (0.35, 0.05)--(0.45,0.15)--(0.55, 0.05);

        \draw[important line] (0.6, -0.02)--(0.6, 0.02);\node[scale=0.6] at(0.6, -0.08){$01(-1)0^\f$};\node[scale=0.6] at(0.6, 0.08){$\frac{4}{27}$};
        \draw[-latex] (0.65, 0.05)--(0.85, 0.15);

           \draw[important line] (0.9, -0.02)--(0.9, 0.02);\node[scale=0.6] at(0.9, -0.08){$01^\f$};\node[scale=0.6] at(0.9, 0.08){$\frac{1}{3}$};
           \draw[-latex] (0.95, 0.15)--(1.15, 0.05);

           \draw[important line] (1.2, -0.02)--(1.2, 0.02);\node[scale=0.6] at(1.2, -0.08){$1^\f$};\node[scale=0.6] at(1.2, 0.08){$1$};

\end{tikzpicture}
\end{center}
\caption{Left: the graph of  the functions $\Pi((i_n), \cdot)$ {for} $(i_n)\in(0^\f, 1^\f)$.   Right:   the {piecewise} monotonicity of $\Pi((i_n), \cdot)$ for $(i_n)\in(0^\f, 1^\f)$, {where the} number  above the sequence $(i_n)$ {denotes}   $t=\Pi((i_n), 1/3)$.}\label{fig:monotonicity-Pi}
\end{figure}

{However,   $\Pi$ is not always monotonic with respect to its second variable (see Figure \ref{fig:monotonicity-Pi}, {left}), which complicates our proofs in many places of the paper. Note that the symbolic space $\set{-1,0,1}^\N$ is symmetric, and    $\Pi((i_n), \cdot)$ and $\Pi((-i_n), \cdot)$ have opposite monotonicity. Moreover, $\Pi(0^\f, \la)\equiv 0$ and $\Pi(1^\f, \la)\equiv 1$. So  it suffices to consider the {piecewise} monotonicity of $\Pi((i_n), \cdot)$ for $(i_n)\in(0^\f, 1^\f)$, {which will be divided} into four subintervals (see Figure \ref{fig:monotonicity-Pi}, {right}):
\[
(0^\f, 001^\f],\quad  (001^\f, 01(-1)0^\f),\quad  [01(-1)0^\f, 01^\f]\quad\textrm{{and}}\quad (01^\f, 1^\f).
\]
\begin{lemma}\label{lem:monotonicity-lambda}
{Let $(i_n)\in (0^\f, 1^\f)$.
\begin{enumerate}[{\rm(i)}]
\item If $(i_n)\in (0^\f, 001^\infty]$,  then $\Pi((i_n), \cdot)$ is  strictly increasing in $(0,1/3]$;
\item If {$ (i_n)\in(001^\infty, 01(-1)0^\infty)$}, then $\Pi((i_n), \cdot)$ is strictly concave in $(0,1/3]$. Furthermore, there exists a unique {$\la_{(i_n)}\in[ {1}/{4}, {1}/{3})$}  such that $\Pi((i_n), \cdot)$ is strictly  increasing in $(0,\la_{(i_n)}]$  and strictly decreasing in $[\la_{(i_n)}, {1}/{3}]$;
\item If {$ (i_n)\in[01(-1)0^\infty, 01^\infty]$}, then $\Pi((i_n), \cdot)$ is strictly increasing in $(0,1/3]$;
\item If $  (i_n)\in (01^\infty, 1^\infty)$, then $\Pi((i_n), \cdot)$ is {strictly} decreasing in $(0,1/3]$.
\end{enumerate}}
\end{lemma}
\begin{proof}
For (i) we take $ (i_n)\in(0^\f, 001^\infty]$. Then   there exists $k\ge 3$ such that $i_1\ldots i_{k}=0^{k-1}1$. By (\ref{eq:derivative}) it follows that
\begin{align*}
  \Pi_2((i_n),\la)  &\ge k\left(\frac{k-1}{k}-\la\right)\la^{k-2}-\sum_{n=k+1}^\f n\left(\frac{n-1}{n}-\la\right)\la^{n-2}\\
  &= \la^{k-2}(k-1-2k\la) > 0,
\end{align*}
where the last inequality follows by {$\la\in(0,1/3)$} and $k\ge 3$. This proves (i).

Next we consider (ii). Take {$(i_n)\in(001^\f,  01(-1)0^\f)$}. Then {$01(-1)^\f\lle (i_n)\prec 01(-1)0^\f$}, and thus $i_1i_2i_3=01(-1)$ and {$i_4i_5\ldots\prec 0^\f$}. By (\ref{eq:derivative}) it follows that
\begin{equation}\label{eq:derivative-1}
\begin{split}
 \Pi_2((i_n),\la)&=2\left(\frac{1}{2}-\la\right)-3\left(\frac{2}{3}-\la\right)\la+\sum_{n=4}^\f i_n n\left(\frac{n-1}{n}-\la\right)\la^{n-2}\\
  &=1-4\la+3\la^2+\sum_{n=4}^\f i_n n\left(\frac{n-1}{n}-\la\right)\la^{n-2},
\end{split}
\end{equation}
and therefore
\begin{equation}\label{eq:secondderivative}
 \frac{\partial\Pi_2((i_n),\la)}{\partial\la}=6\la-4+\sum_{n=4}^\f i_n(n-1)[(n-2)-n\la]\la^{n-3} < 6\la-4< 0.
\end{equation}
So, $\Pi((i_n), \cdot)$ is strictly concave in $(0,1/3]$. Observe by (\ref{eq:derivative-1}) that {$\Pi_2((i_n),\la)\geq 1-4\lambda>0$} for any $\la\in(0,1/4)$, and $\Pi_2((i_n), 1/3)<0$ {since $i_4 i_5 \cdots\prec 0^\infty$}. So by (\ref{eq:secondderivative}) there exists a unique {$\la_{(i_n)}\in[1/4,1/3)$} such that $\Pi((i_n), \cdot)$ is strictly increasing in $(0,\la_{(i_n)}]$ and strictly decreasing in $[\la_{(i_n)}, 1/3]$.

For (iii) we take {$(i_n)\in[01(-1)0^\f, 01^\f]$}. Then $i_1i_2=01$. If $i_3=-1$, then {$i_4i_5\ldots \succeq 0^\f$}. {By} (\ref{eq:derivative}) we obtain that
$
\Pi_2((i_n),\la){\geq}1-4\la+3\la^2> 0$ for any $\la\in(0,1/3)$; and thus
  $\Pi((i_n), \cdot)$ is strictly increasing in {$(0,1/3]$}.
If $i_3\in\set{0,1}$, then by (\ref{eq:derivative}) it follows that
\begin{align*}
\Pi_2((i_n),\la)\geq {\Pi(010(-1)^\infty,\lambda)} = 2\left(\frac{1}{2}-\la\right)-\sum_{n=4}^\f n\left(\frac{n-1}{n}-\la\right)\la^{n-2}
  =1-2\la-3\la^2> 0
  \end{align*}
  for all $\la\in(0,1/3)$.
  This proves (iii).

  Finally  we consider (iv). Take $(i_n)\in(01^\f, 1^\f)$. Then $1(-1)^\f\lle (i_n)\prec 1^\f$, which {gives} $i_1=1$. By (\ref{eq:derivative}) it follows that
\begin{align*}
  \Pi_2((i_n),\la)&=-1+\sum_{n=2}^\f i_n n\left(\frac{n-1}{n}-\la\right)\la^{n-2} < -1+\sum_{n=2}^\f n\left(\frac{n-1}{n}-\la\right)\la^{n-2} =0,
  \end{align*}
  which proves (iv).
\end{proof}
\begin{remark}
  \label{rem:monotonicity-Pi}
From the proof of Lemma \ref{lem:monotonicity-lambda} (ii) it follows that the critical value $\lambda_{(i_n)}$ is {the} unique zero in {$[1/4,1/3)$} of (\ref{eq:derivative-1}). This implies that the map $(i_n)\mapsto \la_{(i_n)}$ is continuous and strictly increasing in $(001^\f, 01(-1)0^\f]=[01(-1)^\f, 01(-1)0^\f]$ with $\la_{01(-1)^\f}=1/4$ and $\la_{01(-1)0^\f}=1/3$.
\end{remark}

{Given $t\in(0,1)\setminus\set{1/3}$, by (\ref{eq:def-Lambda-t}) it follows that
\[\La(t)= \set{\la\in(0,1/3]: \Pi((i_n), \la)=t\textrm{ for some }(i_n)\in(0^\f, 1^\f)}.\]
Note by Lemma  \ref{lem:monotonicity-in} that each $\la\in\La(t)\setminus\set{1/3}$ corresponds to a unique coding $(i_n)$ satisfying $\Pi((i_n),\la)=t$, while for $\la=1/3$ there might be two codings $(i_n)$ satisfying $\Pi((i_n),\la)=t$, one ends with $(-1)^\f$ and the other ends with $1^\f$. {Now we define the map
\begin{equation}\label{eq:Phi-t}
\Phi_t: \La(t) \to(0^\f, 1^\f);\quad {\la\mapsto (i_n) \quad \textrm{with} \quad  \Pi((i_n),\lambda)=t,}
\end{equation}
where for $\lambda=1/3$ we set}
\[\left\{\begin{array}{lll}
\Phi_t(1/3)\textrm{ does not end with }1^\f&\textrm{ if }& t\in(0,1/9]\cup[4/27,1/3),\\
\Phi_t(1/3)\textrm{ does not end with }(-1)^\f&\textrm{ if }& t\in(1/9,4/27)\cup(1/3,1).
\end{array}\right.
\]
So, if $t\in(0,1/9]\cup[4/27,1/3)$, then $\Phi_t(1/3)$ is the \emph{greedy} triadic coding  of $t$; {and} if $t\in(1/9,4/27)\cup(1/3,1)$, then $\Phi_t(1/3)$ is the \emph{lazy} triadic  coding of $t$. {Note that the definition of $\Phi_t(1/3)$ depends on the monotonicity of $\Pi(\Phi_t(1/3),\cdot)$ (see Figure \ref{fig:monotonicity-Pi}, right).}

    In the remaining part of this section we will show that $\Phi_t$ is continuous and piecewise monotonic, which will be vital in our study of $\La(t)$ and its level {sets}.
\begin{lemma}
  \label{lem:continuity-Phi-t}
  For any $t\in(0,1)\setminus\set{1/3}$ the map $\Phi_t$ is continuous in $\Lambda(t)$.
\end{lemma}}
{\begin{proof}
  The proof   is similar to that of \cite[Lemma 2.2]{Wang-Jiang-Kong-Li-2021}. For completeness we sketch its main idea.
  Take $\la_*\in\La(t)$. Suppose on the contrary that $\Phi_t$ is not continuous at $\la_*$. Then there exist  a large $N\in\N$ and a sequence $(\la_k)\subset\La(t)$ such that
  \begin{equation}
    \label{eq:continuity-1}
    \lim_{k\to\f}\la_k=\la_*\quad\textrm{and}\quad {{\rho(\Phi_t(\la_k),\Phi_t(\la_*))}\ge 3^{-N+1}\quad\forall k\ge 1}.
  \end{equation}
  Write $\Phi_t(\la_k)=(i_n^{(k)})$ and $\Phi_t(\la_*)=(i_n^*)$. {By} (\ref{eq:continuity-1}) {it follows that} $i_1^{(k)}\ldots i_N^{(k)}\ne i_1^*\ldots i_N^*$ for all $k\ge 1$. Since $(\set{-1,0,1}^\N, \rho)$ is a compact metric space, there exists a subsequence $(k_j)\subset\N$ such that $\lim_{j\to\f}(i_n^{(k_j)})$ exists, say $(i_n')\in\set{-1,0,1}^\N$. Then
  \begin{equation}
    \label{eq:continuity-1'}
    i_1'\ldots i_N'\ne i_1^*\ldots i_N^*.
  \end{equation}
  Note that
  \begin{equation}
    \label{eq:continuity-2}
    \Pi((i_n^{(k_j)}), \la_{k_j})=t=\Pi((i_n^*), \la_*)\quad\forall j\ge 1.
  \end{equation}
  Letting $j\to\f$ in (\ref{eq:continuity-2}), by using (\ref{eq:continuity-1}) and Lemma \ref{lem:continuity-Pi} it follows that
  \begin{equation}
    \label{eq:continuity-3}
    \Pi((i_n'), \la_*)=t=\Pi((i_n^*), \la_*).
  \end{equation}
  If $\la_*\in(0,1/3)$, then by (\ref{eq:continuity-3}) and Lemma \ref{lem:monotonicity-in} it follows that $(i_n')=(i_n^*)$, leading to a contradiction with (\ref{eq:continuity-1'}).
  If $\la_*=1/3$, then we consider the following two cases.

  Case (I). $t\in(0,1/9]\cup[4/27,1/3)$. Then $(i_n^*)=\Phi_t(1/3)\in(0^\f, 001^\f]\cup[01(-1)0^\f, 01^\f)$. By Lemma \ref{lem:monotonicity-lambda}  it follows that {$\Pi((i_n^*), \cdot)$} is strictly increasing in $(0, 1/3]$. So by (\ref{eq:continuity-2}) and Lemma \ref{lem:monotonicity-in} it follows that
        $        (i_n^{(k_j)}) \lge  (i_n^*)$ for all $j\ge 1,$ and thus
       $(i_n')\lge (i_n^*)$.
        Since $(i_n^*)=\Phi_t(1/3)$ is the greedy triadic coding of $t$, we must have $(i_n')= (i_n^*)$,   contradicting to (\ref{eq:continuity-1'}).

  Case (II). $t\in(1/9,4/27)\cup(1/3,1)$. Then $(i_n^*)=\Phi_t(1/3)\in(001^\f, 01(-1)0^\f)\cup(01^\f, 1^\f)$. Then by Lemma \ref{lem:monotonicity-lambda} there exists $\de>0$ such that {$\Pi((i_n^*), \cdot)$} is strictly decreasing in $(1/3-\de, 1/3]$. So by (\ref{eq:continuity-2}) and Lemma \ref{lem:monotonicity-in} it follows that
     $
        (i_n^{(k_j)}) \lle  (i_n^*)$ for all $j\ge 1,
       $
       and therefore $(i_n')\lle (i_n^*)$.
        Since $(i_n^*)$ is the lazy triadic coding of $t$, we must have $(i_n')= (i_n^*)$, again leading to a contradiction with (\ref{eq:continuity-1'}).

  {Hence}, $\Phi_t$ is continuous at $\la_*$. {Since $\lambda_*\in\Lambda(t)$ {was} arbitrary, $\Phi(t)$ is continuous in $\Lambda(t)$.}
\end{proof}
We will end   this section by showing that $\Phi_t$ is   piecewise monotonic. 

\begin{proposition}[{Key proposition}]
  \label{prop:monotonicity-Phi-t}
  Let $t\in(0,1)\setminus\set{1/3}$.
  \begin{enumerate}
    [{\rm(i)}]
    \item If $t\in(0, {1/9}]$, then  $\Phi_t$ is strictly decreasing in $\La(t)$;

    \item If $t\in(1/9, 4/27)$, then there exist  a unique $\tau{=\tau(t)}\in[1/4,1/3)$ such that $\Phi_t$ is strictly decreasing in $(0,\tau]\cap\La(t)$ and strictly increasing in $[\tau, 1/3]\cap\La(t)$;

    \item If $t\in[{4}/{27}, {1}/{3})$, then $\Phi_t$ is strictly decreasing in $\La(t)$;

    \item If $t\in({1}/{3},1)$, then $\Phi_t$ is strictly increasing in $\La(t)$.
  \end{enumerate}
\end{proposition}

{The} proof of Proposition \ref{prop:monotonicity-Phi-t} {will be split} into several lemmas. First we consider (iii) and (iv).
\begin{lemma}
  \label{lem:mon-Phi-t-4/27-1}
  \begin{enumerate}
    [{\rm(i)}]
   \item If $t\in[{4}/{27}, {1}/{3})$, then $\Phi_t$ is strictly decreasing in $\La(t)$;

    \item If $t\in({1}/{3},1)$, then $\Phi_t$ is strictly increasing in $\La(t)$.
  \end{enumerate}
\end{lemma}
\begin{proof}
{In view of} Lemma \ref{lem:monotonicity-lambda}, the proof of (ii) is very similar to (i). {Here} we only prove (i).
  Take $t\in[4/27, 1/3)$. Then $\Phi_t(1/3)\in[01(-1)0^\f, 01^\f)$. By Lemma \ref{lem:monotonicity-lambda} (iii) it follows that for any $\la\in\La(t)\setminus\set{1/3}$
  \[
  \Pi(\Phi_t(\la),\la)=t=\Pi(\Phi_t(1/3),1/3)>\Pi(\Phi_t(1/3),\la),
  \]
which, together with Lemma \ref{lem:monotonicity-in}, implies   $\Phi_t(\la)\succ\Phi_t(1/3)$. Furthermore, note that $\Phi_t(\la)\lle 01^\f$ for any $\la\in\La(t)\setminus\set{1/3}$, since otherwise $t=\Pi(\Phi_t(\la),\la)\ge \Pi(1(-1)^\f, \la)=1-2\la>1/3$,   a contradiction. So, {$\Phi_t(\la)\in[01(-1)0^\f, 01^\f{]}$} for all $\la\in\La(t)$. By Lemma \ref{lem:monotonicity-lambda} (iii) it follows that for any $\la_1,\la_2\in\La(t)$ with $\la_1<\la_2$   \[
  \Pi(\Phi_t(\la_1),\la_1)=t=\Pi(\Phi_t(\la_2),\la_2)>\Pi(\Phi_t(\la_2),\la_1),
  \]
  which implies $\Phi_t(\la_1)\succ\Phi_t(\la_2)$. This completes the proof.
\end{proof}
In the following we consider $t\in(0,4/27)$. Note by Lemma \ref{lem:monotonicity-lambda} (ii) that $\Pi((i_n), \cdot)$ is not globally  monotonic in $(0,1/3]$ for $(i_n)\in[01(-1)^\f, 01(-1)0^\f)${, which} complicates our proofs of {Proposition \ref{prop:monotonicity-Phi-t} (i) and (ii)}.
{\begin{lemma}
  \label{lem:mon-Phi-t-1}
  Let $t\in(0,4/27)$, and let $\la_\diamond=\la_\diamond(t)$ be the unique root in {$(0,1/3)$} of the equation $\la_\diamond(1-\la_\diamond)^2=t$. Then $\Phi_t$ is strictly decreasing in $(0,\la_\diamond]\cap\La(t)$.
\end{lemma}}
\begin{proof}
The proof is similar to Lemma \ref{lem:mon-Phi-t-4/27-1}.
Since  $t=\la_\diamond(1-\la_\diamond)^2=\Pi(01(-1)0^\f, \la_\diamond)$,
by Lemma \ref{lem:monotonicity-lambda} (iii) it follows that for any $\la\in(0, \la_\diamond]\cap\La(t)$
\[
\Pi(\Phi_t(\la), \la)=t=\Pi(01(-1)0^\f, \la_\diamond)\ge\Pi(01(-1)0^\f, \la),
\]
which gives $\Phi_t(\la)\lge 01(-1)0^\f$. Furthermore, $\Phi_t(\la)\lle 01^\f$ since $t<1/3$. Therefore, $\Phi_t(\la)\in[01(-1)0^\f, 01^\f]$ for any $\la\in(0,\la_\diamond]\cap\La(t)$. By the same argument as in the proof of Lemma \ref{lem:mon-Phi-t-4/27-1} (i) one can   prove that $\Phi_t$ is strictly decreasing in $(0,\la_\diamond]\cap\La(t)$.
\end{proof}
Now we are ready to prove {Proposition \ref{prop:monotonicity-Phi-t} (i)}. Note by Lemma \ref{lem:monotonicity-lambda}  (ii) that for any $(i_n)\in(001^\f, 01(-1)0^\f)$ there exists a unique $\la_{(i_n)}\in[1/4, 1/3)$ such that $\Pi((i_n), \cdot)$ is strictly increasing in $(0,\la_{(i_n)}]$ and strictly decreasing in $[\la_{(i_n)}, 1/3]$.
\begin{lemma}
  \label{lem:mon-phi-t-0-1/9}
  Let $t\in(0,1/9]$. Then $\Phi_t$ is strictly decreasing in $\La(t)$.
\end{lemma}
\begin{proof}
{First we show that $\Phi_t$ is strictly decreasing in $[\sqrt{t}, 1/3]\cap\La(t)$.
    Observe  that  $t=(\sqrt{t})^2=\Pi(001^\f, \sqrt{t})$. By the same argument as {in} the proof of Lemma \ref{lem:mon-Phi-t-1}  we can prove that $\Phi_t(\la)\in(0^\f, 001^\f]$ for any $\la\in[\sqrt{t}, 1/3]\cap\La(t)$, and then by  Lemma \ref{lem:monotonicity-lambda} (i) it follows that $\Phi_t$ is strictly decreasing in $[\sqrt{t}, 1/3]\cap\La(t)$. {Note by Lemma \ref{lem:monotonicity-in} that
    $\Pi(01(-1)0^\f, \la_\diamond)=t=\Pi(001^\f, \sqrt{t})<\Pi(01(-1)0^\f, \sqrt{t}).$
    Then by Lemma \ref{lem:monotonicity-lambda} (iii) it follows that $\la_\diamond<\sqrt{t}$. So,}
 by Lemma \ref{lem:mon-Phi-t-1} it suffices to prove that $\Phi_t$ is strictly decreasing in $(\la_\diamond, \sqrt{t})\cap\La(t)$.}
       Note that $\Pi(01(-1)0^\f,\la_\diamond)=t=\Pi(001^\f,\sqrt{t})$. Then by Lemma \ref{lem:monotonicity-lambda} it follows that for any $\la\in(\la_\diamond, \sqrt{t})\cap\La(t)$
  \[
  \Pi(001^\f, \la)<\Pi(001^\f, \sqrt{t})=t=\Pi(\Phi_t(\la),\la)=\Pi(01(-1)0^\f, \la_\diamond)<\Pi(01(-1)0^\f, \la),
  \]
  which implies
  \begin{equation}\label{eq:dec14-2}
  \Phi_t(\la)\in(001^\f, 01(-1)0^\f).
  \end{equation}

  Now we claim that $\la<\la_{\Phi_t(\la)}$ for all $\la\in (\la_\diamond, \sqrt{t})\cap\La(t)$.
  Suppose on the contrary that $\la\ge \la_{\Phi_t(\la)}$. Then by the definition of $\la_{\Phi_t(\la)}$ it follows that for any $\la\in(\la_\diamond, \sqrt{t})\cap\La(t)$
  \[
  \Pi(\Phi_t(\la),\sqrt{t})<\Pi(\Phi_t(\la), \la)=t=\Pi(001^\f, \sqrt{t}),
  \]
  which yields $\Phi_t(\la)\prec 001^\f$, leading to a contradiction with (\ref{eq:dec14-2}). This proves the claim.

  {Therefore,} for any $\la_1,\la_2\in(\la_\diamond, \sqrt{t})\cap\La(t)$ with $\la_1<\la_2$
  \[
  \Pi(\Phi_t(\la_1),\la_1)=t=\Pi(\Phi_t(\la_2),\la_2)>\Pi(\Phi_t(\la_2), \la_1),
  \]
  {where the inequality follows by $\lambda_1<\lambda_2<\lambda_{\Phi_t({\lambda_2})}$.
  This} gives $\Phi_t(\la_1)\succ\Phi_t(\la_2)$, completing the proof.
\end{proof}
In the following it {remains} to prove {Proposition \ref{prop:monotonicity-Phi-t} (ii)}. First we consider $t\in(1/9,1/8]$.
\begin{lemma}
  \label{lem:mon-phi-t-1/9-1/8}
  Let $t\in(1/9,1/8]$. Then $\Phi_t$ is strictly decreasing in $(0,1/4]\cap\La(t)$ and strictly increasing in $[1/4, 1/3]\cap\La(t)$.
\end{lemma}
\begin{proof}
{Note that $\Pi(01(-1)0^\f, \la_\diamond)=t\le  {1}/{8}< {9}/{64}=\Pi(01(-1)0^\f,  {1}/{4}).$ Then by Lemma \ref{lem:monotonicity-lambda} (iii) it follows that $\la_\diamond=\la_\diamond(t)<1/4$ for all $t\in(1/9, 1/8]$. So,}
  by Lemma \ref{lem:mon-Phi-t-1}   it suffices to prove that $\Phi_t$ is strictly decreasing in {$(\la_\diamond, 1/4]\cap\Lambda(t)$} and strictly increasing in {$[1/4,1/3]\cap\Lambda(t)$}. Note that $\Pi(01(-1)0^\f, \la_\diamond)=t>1/9=\Pi(001^\f, 1/3)$. Then by Lemma \ref{lem:monotonicity-lambda} it follows that
  for any $\la\in(\la_\diamond, 1/3)\cap\La(t)$
  \[
 \Pi(001^\f, \la)\le\Pi(001^\f, 1/3) <\Pi(\Phi_t(\la), \la)=t=\Pi(01(-1)0^\f, \la_\diamond)<\Pi(01(-1)0^\f,\la).
  \]
  This implies that
  \begin{equation}\label{eq:dec15-1}
  \Phi_t(\la)\in(001^\f, 01(-1)0^\f)=[01(-1)^\f, 01(-1)0^\f).
   \end{equation}
   So, by Lemma \ref{lem:monotonicity-lambda} (ii) it follows that   for any $\la_1, \la_2\in\La(t)\cap(\la_\diamond, 1/4]$ with $\la_1<\la_2$
  \[
  \Pi(\Phi_t(\la_1), \la_1)=t=\Pi(\Phi_t(\la_2), \la_2)>\Pi(\Phi_t(\la_2),\la_1),
  \]
  which gives $\Phi_t(\la_1)\succ \Phi_t(\la_2)$.

  Next we show that $\Phi_t$ is strictly increasing in $[1/4, 1/3]\cap\La(t)$. We claim that {$\la> \la_{\Phi_t(\la)}$} in this case. Suppose on the contrary that $\la\le\la_{\Phi_t(\la)}$ for {some} $\la\in(1/4,1/3)\cap\La(t)$. Then
  \[
  \Pi(\Phi_t(\la), 1/4)<\Pi(\Phi_t(\la), \la)=t\le 1/8=\Pi(01(-1)^\f, 1/4),
  \]
  which implies $\Phi_t(\la)\prec 01(-1)^\f$, leading to a contradiction with (\ref{eq:dec15-1}). This proves the claim, and then it follows that for any $\la_1, \la_2\in(1/4,1/3)\cap\La(t)$ with $\la_1<\la_2$
  \[
  \Pi(\Phi_t(\la_2),\la_2)=t=\Pi(\Phi_t(\la_1),\la_1)>\Pi(\Phi_t(\la_1),\la_2),
  \]
  where the inequality follows {by} {$\lambda_2>\lambda_1> \lambda_{\Phi_t(\lambda_1)}$}. This gives $\Phi_t(\la_2)\succ\Phi_t(\la_1)$.
\end{proof}
To prove {Proposition \ref{prop:monotonicity-Phi-t} (ii)} for $t\in(1/8, 4/27)$ we need  two lemmas on the critical value $\la_{(i_n)}$.
\begin{lemma}\label{lem:mon-phi-t-2}
The map $\phi: (i_n)\mapsto\Pi((i_n), \la_{(i_n)})$ is strictly increasing and continuous in  the interval $[01(-1)^\f, 01(-1)0^\f)$.
\end{lemma}
\begin{proof}
Note by Remark \ref{rem:monotonicity-Pi} that $(i_n)\mapsto \la_{(i_n)}$ is continuous. Then the continuity of $\phi$ follows by Lemma \ref{lem:continuity-Pi}. And the monotonicity of $\phi$ follows by Lemma \ref{lem:monotonicity-in} that for any $(i_n), (j_n)\in[01(-1)^\f, 01(-1)0^\f)$ with $(i_n)\prec (j_n)$
\[
\Pi((i_n), \la_{(i_n)})<\Pi((j_n), \la_{(i_n)})\le \Pi((j_n), \la_{(j_n)}),
\]
where the last inequality follows by Lemma \ref{lem:monotonicity-lambda} (ii) that $\Pi((j_n),\cdot)$ attains its maximum value at {$\la_{(j_n)}\in[1/4,1/3)$}.
\end{proof}

\begin{lemma}
  \label{lem:mon-phi-t-3}
  Let $t\in(1/8, 4/27)$. Then there exists a unique {$\tau=\tau(t)\in(1/4,1/3)\cap\Lambda(t)$} such that $\la_{\Phi_t(\tau)}=\tau$. Furthermore, the map $t\mapsto\tau(t)$ is strictly increasing.
\end{lemma}
\begin{proof}
  Note by Remark \ref{rem:monotonicity-Pi} that $\la_{01(-1)^\f}=1/4$ and $\la_{01(-1)0^\f}=1/3$.   Then
  \[
  \phi(01(-1)^\f)=\Pi\left(01(-1)^\f,\frac{1}{4}\right)=\frac{1}{8},\quad \phi(01(-1)0^\f)=\Pi\left(01(-1)0^\f, \frac{1}{3}\right)=\frac{4}{27}.
  \]
   {So, by} Lemma \ref{lem:mon-phi-t-2} it follows that $\phi$ bijectively maps {$(01(-1)^\f, 01(-1)0^\f)$ to $(1/8,4/27)$}.  Since $t\in(1/8, 4/27)$,  by Lemma \ref{lem:mon-phi-t-2} there exists a unique $(j_n)\in(01(-1)^\f, 01(-1)0^\f)$ such that
  \begin{equation}
    \label{eq:dec14-1}
    t=\Pi((j_n), \la_{(j_n)}).
  \end{equation}
  Observe by Remark \ref{rem:monotonicity-Pi} that the map $(i_n)\mapsto \la_{(i_n)}$ is  strictly increasing, $\la_{01(-1)^\f}=1/4$ {and} $\la_{01(-1)0^\f}=1/3$. Then by (\ref{eq:dec14-1}) it follows that $\tau:=\la_{(j_n)}\in(1/4, 1/3)\cap\La(t)$ and $(j_n)=\Phi_t(\tau)$.     Thus, $\la_{\Phi_t(\tau)}=\la_{(j_n)}=\tau$.

  Next we prove that the map $t\mapsto \tau(t)$ is strictly increasing.
  {Take {$t_1,t_2\in(1/8,4/27)$} with $t_1<t_2$. Write $\tau_i=\tau(t_i)$ for $i=1, 2$. Since
  $
  \Pi(\Phi_{t_1}(\tau_1),\tau_1)=t_1<t_2=\Pi(\Phi_{t_2}(\tau_2),\tau_2),
 $
  by Lemma \ref{lem:mon-phi-t-2} we have $\Phi_{t_1}(\tau_1)\prec\Phi_{t_2}(\tau_2)$. Hence,   by Remark \ref{rem:monotonicity-Pi} it follows that
  \[\tau_1=\la_{\Phi_{t_1}(\tau_1)}<\la_{\Phi_{t_2}(\tau_2)}=\tau_2\]
  as desired.}
\end{proof}
\begin{lemma}
  \label{lem:mon-phi-t-1/8-4/27}
  Let $t\in(1/8,4/27)$. Then $\Phi_t$ is strictly decreasing in $(0,\tau]\cap\La(t)$ and strictly increasing in $[\tau, 1/3]\cap\La(t)$, where $\tau=\tau(t)$ is defined as in Lemma \ref{lem:mon-phi-t-3}.
\end{lemma}
\begin{proof}
{Note by the proof of Lemma \ref{lem:mon-phi-t-3} that $\Phi_t(\tau)\in(01(-1)^\f, 01(-1)0^\f)$ for any $t\in(1/8,4/27)$. Then by Lemma \ref{lem:monotonicity-in} it follows that
$
\Pi(01(-1)0^\f, \la_\diamond)=t=\Pi(\Phi_t(\tau),\tau)<\Pi(01(-1)0^\f, \tau),
$
which implies $\la_\diamond<\tau$ by Lemma \ref{lem:monotonicity-lambda} (iii). So,}
  by Lemma \ref{lem:mon-Phi-t-1}   it suffices to prove that $\Phi_t$ is strictly decreasing in $(\la_\diamond, \tau]\cap\La(t)$ and strictly increasing in $[\tau, 1/3]\cap\La(t)$.
  First we claim that $\la\le\la_{\Phi_t(\la)}$ for all $\la\in(\la_\diamond,\tau]\cap\La(t)$. Note that for $\la\in(\la_\diamond,\tau]\cap\La(t)$
  \[
  \Pi(\Phi_t(\la), \la)=t=\Pi(\Phi_t(\tau),\tau)\ge\Pi(\Phi_t(\tau),\la),
  \]
  which gives $\Phi_t(\la)\lge\Phi_t(\tau)$. Since the map $(i_n)\mapsto\la_{(i_n)}$ is strictly increasing,   by Lemma \ref{lem:mon-phi-t-3} we conclude that $\la_{\Phi_t(\la)}\ge\la_{\Phi_t(\tau)}=\tau\ge\la$, proving the claim. Therefore, for any $\la_1,\la_2\in(\la_\diamond,\tau]\cap\La(t)$ with $\la_1<\la_2$
  \[
  \Pi(\Phi_t(\la_1), \la_1)=t=\Pi(\Phi_t(\la_2),\la_2)>\Pi(\Phi_t(\la_2),\la_1),
  \]
  which yields $\Phi_t(\la_1)\succ\Phi_t(\la_2)$.

  Next we prove that $\Phi_t$ is strictly increasing in $[\tau, 1/3]\cap\La(t)$. We claim that $\la\ge \la_{\Phi_t(\la)}$  in this case. Suppose on the contrary that $\la<\la_{\Phi_t(\la)}$ for some $\la\in[\tau, 1/3]\cap\La(t)$. Then
  \[
  \Pi(\Phi_t(\tau),\tau)=t=\Pi(\Phi_t(\la),\la)\ge\Pi(\Phi_t(\la),\tau),
  \]
  which implies that $\Phi_t(\tau)\lge\Phi_t(\la)$. Since the map $(i_n)\mapsto\la_{(i_n)}$ is strictly increasing, by Lemma \ref{lem:mon-phi-t-3} it follows that $\la_{\Phi_t(\la)}\le\la_{\Phi_t(\tau)}=\tau\le\la$, leading to a contradiction. This proves the claim,  {and then} it follows that for any $\la_1,\la_2\in[\tau, 1/3]\cap\La(t)$ with $\la_1<\la_2$
  \[
  \Pi(\Phi_t(\la_2),\la_2)=t=\Pi(\Phi_t(\la_1),\la_1)>\Pi(\Phi_t(\la_1),\la_2),
  \]
  which gives $\Phi_t(\la_2)\succ\Phi_t(\la_1)$. This completes the proof.
\end{proof}
\begin{proof}[Proof of Proposition \ref{prop:monotonicity-Phi-t}]
The proposition follows by {Lemmas \ref{lem:mon-Phi-t-4/27-1}, \ref{lem:mon-phi-t-0-1/9}, \ref{lem:mon-phi-t-1/9-1/8} and \ref{lem:mon-phi-t-1/8-4/27}}.
\end{proof}
{By the proof of Proposition \ref{prop:monotonicity-Phi-t} (i) and (ii) it follows that
\begin{equation}\label{eq:la-diamond-tau}
t<\la_\diamond(t)<\sqrt{t}\le \frac{1}{3}~\textrm{ if }0<t\le\frac{1}{9},\qquad t<\la_\diamond(t)<\tau(t)<\frac{1}{3}~\textrm{ if } \frac{1}{9}<t< \frac{4}{27}.
\end{equation}
Here $\tau(t)=1/4$ for any $t\in(1/9, 1/8]$, and $\tau(t)$ is defined as in Lemma \ref{lem:mon-phi-t-3} for $t\in(1/8,4/27)$.}

\section{Topology and Hausdorff dimension of $\La(t)$}\label{sec:topology-La(t)}
In this section we will show that $\La(t)$ is a topological Cantor set {(Proposition \ref{prop:topology-La(t)}) and it has} full Hausdorff dimension ({Corollary \ref{cor:full-dimension}}).

\subsection{Topology of $\La(t)$} Recall that a topological Cantor set in $\R$ is a non-empty perfect set with no interior points.
\begin{proposition}
  \label{prop:topology-La(t)}
  For any $t\in(0,1)\setminus\set{1/3}$ the set $\La(t)$ is a topological Cantor set with
  \[
  \min\La(t)=\min\set{t,\frac{1-t}{2}} \quad\textrm{and}\quad \max\La(t)=\frac{1}{3}.
  \]
\end{proposition}
}
First we determine the extreme points of $\La(t)$.
\begin{lemma}\label{lem:min-max-La(t)}
For any   $t\in(0,1)\setminus \{1/3\}$  we have
$\min\La(t)=\min\{t,\frac{1-t}{2}\}$ and $\max(\La(t))=1/3$.
\end{lemma}
\begin{proof}
Note that $E_\la=[-1,1]$ for $\la=1/3$. This implies that $\max\La(t)=1/3$ for any $t\in[-1, 1]$. For the minimum {value} of $\La(t)$ we consider     two cases: (I) $t\in(0,1/3)$; (II) $t\in(1/3, 1)$.

(I) $t\in(0,1/3)$. Then $\Phi_t(\la)\in(0^\f, 01^\f]$ for any $\la\in(0,1/3]\cap\La(t)$. By Lemmas \ref{lem:mon-Phi-t-4/27-1} (i) and \ref{lem:mon-Phi-t-1}   it follows that the minimum value  $\la_*=\min\La(t)$ satisfies
$
\Pi(01^\f, \la_*)=t.
$
So, $\min\La(t)=\la_*=t$.

(II) $t\in(1/3, 1)$. Then $\Phi_t(\la)\in[1(-1)^\f, 1^\f]$ for any $\la\in\La(t)$. Then by Lemma \ref{lem:mon-Phi-t-4/27-1} (ii) it follows that $\la_*=\min\La(t)$ satisfies $\Pi(1(-1)^\f, \la_*)=t$, which implies that
$t=1-2\la_*$. Hence, $\min\La(t)=\la_*=\frac{1-t}{2}$.
\end{proof}

Next we show that $\La(t)$ is a topological Cantor set.
\begin{lemma}
  \label{lem:perfect-La(t)}
  For any $t\in(0,1)\setminus\set{1/3}$,  $\La(t)$ is a non-empty perfect set.
\end{lemma}
\begin{proof}
By Lemma \ref{lem:min-max-La(t)} it suffices to prove that $\La(t)$ is closed and has no isolated points. First we prove the closeness of $\La(t)$. Let $(\la_j)\subset\La(t)$ with $\lim_{j\to\f}\la_j=\la_0$.
Suppose on the contrary that $\la_0\notin\La(t)$, i.e., $t\notin E_{\la_0}$. Since $E_{\la_0}$ is compact, we have
$
dist(t, E_{\la_0}):=\inf\set{|x-t|: x\in E_{\la_0}}>0.
$
 Note that $E_{\la_0}=\bigcap_{n=1}^\f E_{\la_0}(n)$,
where {for $\la\in(0,1/3]$}
\begin{equation}\label{eq:E_n}
E_\la(n):=\bigcup_{i_1\ldots i_n\in\set{-1,0,1}^n}g_{i_1}\circ\cdots\circ g_{i_n}([-1,1])
\end{equation}
{with} $g_i(x)=\la x+i(1-\la)$ for $i=-1, 0,1$. Since $E_{\la_0}(n)\supset E_{\la_0}(n+1)$ for all $n\ge 1$,  there exists  $N\in\N$ such that
\begin{equation}
  \label{eq:compactness-2}
  d_0:=dist(t, E_{\la_0}(N))>0.
\end{equation}
Note that $\lim_{j\to\f}\la_j=\la_0$, and by (\ref{eq:E_n}) each $E_{\la}(N)$ is the union of $3^N$ pairwise disjoint intervals of equal length $2\la^N$. Then    $E_{\la_j}(N)$ converges to $E_{\la_0}(N)$ under the Hausdorff metric $d_H$ as $j\to\f$. So there exists a large $j\in\N$ such that
\begin{equation*}
  \label{eq:compactness-3}
 dist(t, E_{\la_0}(N))\le  d_H(E_{\la_j}(N), E_{\la_0}(N))\le \frac{d_0}{2},
\end{equation*}
where the first inequality follows by $t\in E_{\la_j}\subset E_{\la_j}(N)$. This leads to a contradiction with (\ref{eq:compactness-2}). So, $\la_0\in \La(t)$, and thus $\La(t)$ is closed.

Next we prove that $\La(t)$ has no isolated points. Take $\la\in\La(t)$ and let $(i_n)=\Phi_t(\la)$. Then there exists a subsequence $(n_j)\subset\N$ such that $i_{n_j}\in\set{-1,0}$ for all $j\ge 1$, or $i_{n_j}\in\set{0,1}$ for all $j\ge 1$. Without lose of generality we assume { that $n_1$ is sufficiently large, and}  $i_{n_j}\in\set{-1,0}$ for all $j\ge 1$. For any $j\ge 1$ we define {$\la_j\in(0,1/3]\cap (\lambda-\delta, \lambda+\delta)$ for some small $\delta >0$,} such that
      $
        \Phi_t(\la_j)=i_1\ldots i_{n_j-1} (i_{n_j}+1)0^\f.
     $
        Then it is clear that $(\la_j)\subset\La(t)$ and $\la_j\to\la$ as $j\to\f$. Hence, $\la$ is not isolated in $\La(t)$.
\end{proof}

\begin{proof}
  [Proof of Proposition \ref{prop:topology-La(t)}]
  By Lemmas \ref{lem:min-max-La(t)} and \ref{lem:perfect-La(t)} it suffices to prove that $\La(t)$ has no interior points. Take $\la_1, \la_2\in\La(t)$ with $\la_1<\la_2$. It suffices to prove   $(\la_1,\la_2)\setminus\La(t)\ne\emptyset$. By Proposition  \ref{prop:monotonicity-Phi-t} there exists $\hat\la\in(\la_1,\la_2)$ such that  $\Phi_t$ is strictly monotonic in $(\hat\la, \la_2)$. Without lose of generality we may assume that $\Phi_t$ is strictly increasing in $(\hat\la, \la_2)$. Write $(i_n)=\Phi_t(\hat\la)$ and $(j_n)=\Phi_t(\la_2)$. Then $(i_n)\prec (j_n)$. So there exists $N\in\N$ such that $i_1\ldots i_{N-1}=j_1\ldots j_{N-1}$ and $i_N<j_N$. {{Suppose} $j_N-i_N=1$, since otherwise we can choose a larger $\hat\la$.}  Let
  \[
  \c=i_1\ldots i_N 1^\f\quad\textrm{and}\quad \d=j_1\ldots j_N(-1)^\f.
  \]
  Then by Lemma \ref{lem:monotonicity-in} it follows that
  \begin{equation}
    \label{eq:no-interiorpoint-La(t)}
    \Pi(\c, \la_2)<\Pi(\d,\la_2)\le \Pi((j_n), \la_2)=t=\Pi((i_n), \hat\la)\le\Pi(\c,\hat\la)<\Pi(\d,\hat\la).
  \end{equation}
  Denote the open interval by $I_\la:=(\Pi(\c,\la), \Pi(\d,\la))$. Observe that the map $\la\mapsto cl(I_\la)$ is continuous with respect to the Hausdorff metric $d_H$. Then   (\ref{eq:no-interiorpoint-La(t)}) implies that $t\in I_\la$ for some $\la\in(\hat\la,\la_2)\subset(\la_1,\la_2)$. Since $I_\la\cap E_\la=\emptyset$, we have $t\notin E_\la$, i.e., $\la\notin\La(t)$. So, $(\la_1,\la_2)\setminus\La(t)\ne\emptyset$, completing the proof.
\end{proof}

\subsection{Hausdorff dimension of $\Lambda(t)$}\label{sec:full-dimension}
Now we turn to prove that  $\La(t)$ has full Hausdorff dimension, which can  be deduced from  the following {result}.

\begin{proposition}\label{prop:local-dim-1/3}
For any $t\in(0,1)\setminus \{1/3\}$  we have
\[{
\lim_{\de\to 0^+}\dim_H\big(\La(t)\cap (1/3-\delta,1/3+\delta )\big)=1.}
\]
\end{proposition}

Our strategy to prove Proposition \ref{prop:local-dim-1/3} is to construct a sequence of  subsets of $\La(t)\cap(1/3-\delta,1/3+\delta)$ whose Hausdorff dimension can be arbitrarily close to one.
In view of $\Phi_t(1/3)$ defined in (\ref{eq:Phi-t}), we consider   two cases:  {$t\in(0,1/9]\cup[4/27, 1/3)$} and {$t\in(1/9,4/27)\cup(1/3,1)$}.

First we consider  {$t\in(0,1/9]\cup[4/27, 1/3)$}. Then {$\Phi_t(1/3)\in(0^\f, 001^\f]\cup[01(-1)0^\f, 01^\f)$}, and $\Phi_t(1/3)$ does not end with $1^\f$. So there exists a {subsequence} {$\{n_j\}\subset\mathbb{N}_{\geq4}$}  such that $x_{n_j}<1$  for {all} $j\ge 1$. Note by Proposition  \ref{prop:monotonicity-Phi-t} (i) and (iii)  that $\Phi_t$ is strictly  decreasing in $\La(t)$.
Then for each $j\ge 1$ there exists a unique $\eta_j\in(0,1/3)$ such that
\[
t=\Pi(x_1\ldots x_{n_j}^+1^\f, \eta_j).
\]
Accordingly, for $k\ge 1$ we {let}
\begin{equation*}
\Delta^+_{k,j}(t):=\{\la\in\La(t): \Phi_t(\la)=x_1\dots x_{n_j}^+d_1d_2\dots \textrm{ with }d_{ik}\ne -1 ~\forall i\ge 1\}.
\end{equation*}Then by Proposition  \ref{prop:monotonicity-Phi-t} {(i) and (iii)}   it follows that
\begin{equation}\label{eq:subset-Gamma-kj(t)+}
\Delta^+_{k,j}(t)\subset\La(t)\cap(\eta_j,1/3)~\forall k\ge 1;\quad \textrm{ and }\quad \eta_j\nearrow 1/3~\textrm{ {as} } j\to\f.
\end{equation}

\begin{lemma}\label{lem:lip-2}
Let {$t\in(0,1/9]\cup[4/27, 1/3)$}. Then  for any $k, j\in\N$ there exists {a constant} $C>0$  such that
\[
|\pi_{\eta_j}(\Phi_t(\la_1))-\pi_{\eta_j}(\Phi_t(\la_2))|\leq C|\la_1-\la_2|
\]
for any $\la_1,\la_2\in \Delta^+_{k,j}(t)$.
\end{lemma}
\begin{proof}
Let $\la_1,\la_2\in \Delta^+_{k,j}(t)$ with $\la_1<\la_2$. Write $(a_n)=\Phi_t(\la_1)$ and $(b_n)=\Phi_t(\la_2)$. Since $\Phi_t$ is strictly decreasing in $[\eta_j, 1/3]\cap\La(t)$, we have $(a_n)\succ(b_n)$; and then there exists a $N>n_j$ such that $a_1\cdots a_{N-1}=b_1\cdots b_{N-1}$ and $a_N>b_N$.  {{So,}
\begin{equation}\label{eq:lip-2}
\begin{split}
|\pi_{\eta_j}(\Phi_t(\lambda_1))-\pi_{\eta_j}(\Phi_t(\lambda_2))|  & =\left|(1-\eta_j) \left(\sum_{n=N}^\infty  a_n\eta_j^{n-1}-\sum_{n=N}^\infty  b_n\eta_j^{n-1}\right)\right| \\
  & \leq (1-\eta_j)\sum_{n=N}^{\infty}2\eta_j^{n-1} = 2\eta_j^{N-1}.
\end{split}
\end{equation}}

On the other hand,  by Lemma \ref{lem:monotonicity-in} it follows that
\[
(1-\la_1)\sum_{n=1}^\infty a_n\la_1^{n-1}=t=(1-\la_2)\sum_{n=1}^\infty b_n\la_2^{n-1}\le (1-\la_2)\left(\sum_{n=1}^N a_n\la_2^{n-1}- \sum_{n=N+1}^\infty\la_2^{n-1}\right),
\]
which  implies
\begin{equation}\label{eq:la_1,la_2}
\begin{split}
(1-\la_1)\sum_{n=N+1}^\infty (a_n+1)\la_1^{n-1}&\le (1-\la_2)\left(\sum_{n=1}^N a_n\la_2^{n-1}-\sum_{n=N+1}^\infty\la_2^{n-1}\right)\\
&\qquad-(1-\la_1)\left(\sum_{n=1}^N a_n\la_1^{n-1}-\sum_{n=N+1}^\infty\la_1^{n-1}\right)\\
&=\Pi(a_1\ldots a_N(-1)^\f, \la_2)-\Pi(a_1\ldots a_N(-1)^\f, \la_1)\\
&\leq 2(\la_2-\la_1),
\end{split}
\end{equation}
where the last inequality follows by (\ref{eq:derivative}) that
$
\left|\Pi_2((i_n),\la)\right|\le  2
$ for any $(i_n)\in(0^\f, 1^\f)$ {and} $\la\in(0,1/3]$.
{Note by the definition of $\Delta_{k,j}^+(t)$ that} $a_{N+1} a_{N+2}\cdots\lge ((-1)^{k-1}0)^\f$. {Then} $(a_{N+1}+1)(a_{N+2}+1)\cdots\lge (0^{k-1}1)^\f$. {So,} by  (\ref{eq:la_1,la_2}) and using $\la_1\in(\eta_j,1/3)$ it follows that
\[
2(\la_2-\la_1)>(1-\la_1)\la_1^{N+k-1}>\frac{2}{3}\eta_j^{N+k-1}.
\]
 This, together with  (\ref{eq:lip-2}), implies that
\[
{|\pi_{\eta_j}(\Phi_t(\lambda_1))-\pi_{\eta_j}(\Phi_t(\lambda_2))|}\leq \frac{6}{\eta_j^k} |\la_1-\la_2|
\]
as required.
\end{proof}

Next we consider {$t\in(1/9,4/27)\cup(1/3,1)$}. Then $\Phi_t(1/3)\in(001^\f, 01(-1)0^\f)\cup(01^\f, 1^\f)$, {and it} does not end with $(-1)^\f$. So there exists a {subsequence} $\{n_j\}\subset\mathbb{N}$  such that $x_{n_j}>-1$  for any $j\ge 1$.  By Proposition  \ref{prop:monotonicity-Phi-t} (ii) and (iv)  there exists $\de>0$ such that $\Phi_t$ is strictly  increasing in $\La(t)\cap(1/3-\de, 1/3]$.
So, without loss of generality, by deleting the first {finitely many} terms of $\set{n_j}$, we can assume   that for any $j\ge 1$  the equation
\[
t=\Pi(x_1\ldots x_{n_j}^-(-1)^\f, \gamma_j)
\]
determines a unique  $\gamma_j\in(1/3-\de,1/3)$.
Note by Lemma \ref{lem:monotonicity-lambda} (ii) that if $x_1\ldots x_{n_j}^-(-1)^\f\in(001^\f, 01(-1)0^\f)$, the above equation  {may determine} two different $\gamma_j$s in $(0,1/3)$, but  only one is in $(1/3-\de, 1/3)$.
Accordingly, for $k\ge 1$ we {set}
\begin{equation*}
\Delta^-_{k,j}(t):=\{\la\in\La(t)\cap(1/3-\de, 1/3): \Phi_t(\la)=x_1\dots x_{n_j}^-d_1d_2\dots \textrm{ with }d_{ik}\ne -1 ~\forall i\ge 1\}.
\end{equation*}Then by Proposition \ref{prop:monotonicity-Phi-t} {(ii) and (iv)}   it follows that
\begin{equation}\label{eq:subset-Gamma-kj(t)-}
\Delta_{k,j}^-(t)\subset\La(t)\cap(\gamma_j,1/3)~\forall k\ge 1;\quad \textrm{ and }\quad \gamma_j\nearrow 1/3~\textrm{ as }j\to\f.
\end{equation}

\begin{lemma}\label{lem:lip-1}
Let {$t\in(1/9,4/27)\cup(1/3,1)$}. Then for any $k,j\in\mathbb{N}$  there exists  a constant $C>0$  such that
\[{
|\pi_{\gamma_j}(\Phi_t(\la_1))-\pi_{\gamma_j}(\Phi_t(\la_2))|\leq C|\la_1-\la_2|}
\]
for any $\la_1,\la_2\in \Delta^-_{k,j}(t)$.
\end{lemma}

\begin{proof}
Note that
  $\Phi_t$ is strictly increasing in $[\gamma_j, 1/3]\cap\La(t)$.
Then by a similar argument as in the proof of Lemma \ref{lem:lip-2} one can {verify} that
$
|\pi_{\gamma_j}(\Phi_t(\la_1))-\pi_{\gamma_j}(\Phi_t(\la_2))|\leq \frac{6}{\gamma_j^k} |\la_1-\la_2|$ for any $\la_1,\la_2\in\Delta_{k,j}^-(t).
$
\end{proof}

\begin{proof}
  [Proof of Proposition \ref{prop:local-dim-1/3}]
  Take $t\in(0,1)\setminus\set{1/3}$. Since the proof for {$t\in(1/9,4/27)\cup(1/3,1)$} is similar, we  {only consider}  $t\in(0,1/9]\cup[4/27,1/3)$. {By} (\ref{eq:subset-Gamma-kj(t)+}) and Lemma \ref{lem:lip-2} it follows that
  \begin{equation}
    \label{eq:localdim-1/3-1}
    \dim_H(\La(t)\cap(\eta_j,1/3))\ge \dim_H\Delta_{k,j}^+(t)\ge \dim_H\pi_{\eta_j}(\Phi_t(\Delta_{k,j}^+(t)))
  \end{equation}
  for {all $k, j\in\N$}. Note by the definition of $\Delta_{k,j}^+(t)$ that $\Phi_t(\Delta_{k,j}^+(t))$ consists of all sequences $(i_n)\in\set{-1,0,1}^\N$ with a prefix $i_1\ldots i_{n_j}=x_1\ldots x_{n_j}^+$ and  $i_{n_j+mk}\ne -1$ for all $m\ge 1$. So, by using $\eta_j\in(0,1/3)$ and (\ref{eq:localdim-1/3-1}) we obtain that
  \begin{equation}
    \label{eq:localdim-1/3-2}
    \dim_H(\La(t)\cap(\eta_j,1/3))\ge\frac{(k-1)\log 3+\log 2}{-k\log\eta_j}\to \frac{\log3}{-\log\eta_j} {\quad \textrm{as } k\to \infty.}
  \end{equation}
  Hence, by (\ref{eq:subset-Gamma-kj(t)+}) and (\ref{eq:localdim-1/3-2}) it follows that
  \[
 \lim_{\de\to 0^+} \dim_H(\La(t)\cap(1/3-\de, 1/3+\de))\ge\lim_{j\to\f}\dim_H(\La(t)\cap(\eta_j, 1/3))\ge\lim_{j\to\f}\frac{\log 3}{-\log\eta_j}=1.
  \]
  The reverse inequality is {obvious,} since $\dim_H(\La(t)\cap(1/3-\de, 1/3+\de))\le 1$ for any $\de>0$. 
\end{proof}
{The following corollary follows directly from Proposition \ref{prop:local-dim-1/3}.}
\begin{corollary}
 \label{cor:full-dimension}
 For any $t\in(0,1)\setminus\set{1/3}$ we have $\dim_H\La(t)=1$.
\end{corollary}

\section{Proof of Theorem \ref{thm: La(t)}}\label{sec:local-dimension}
In this section we will determine the local dimension of $\La(t)$, and prove Theorem \ref{thm: La(t)}.
\begin{proposition}
  \label{th:local-dim-lambda}
  Let $t\in(0,1)\setminus\set{1/3}$. Then for any {$\la\in\La(t)$}
  \begin{equation}\label{eq:localdim-lambda}
  \lim_{\de\to 0^+}\dim_H\big(\La(t)\cap(\la-\de, \la+\de)\big)=\frac{\log 3}{-\log\la}.
  \end{equation}
\end{proposition}
  Note by Proposition \ref{prop:local-dim-1/3}   that (\ref{eq:localdim-lambda}) holds for $\la=1/3$. So, in the following we only need to prove (\ref{eq:localdim-lambda}) for $\la\in\La(t)\setminus\set{1/3}$, which will be split into   two subsections.

 \subsection{{A lower bound {on} the local dimension of $\Lambda(t)$}} Take {$\la_*\in\La(t)\setminus\{1/3\}$} and let $(x_i)=\Phi_t(\la_*)$ be the unique  coding of $t$ defined as in (\ref{eq:Phi-t}), i.e., $t=\Pi((x_i), \la_*)$. We will {prove} in this {subsection} that
{\[\lim_{\de\to 0^+}\dim_H\big(\La(t)\cap(\la_*-\de, \la_*+\de)\big)\ge \frac{\log 3}{-\log\la_*}.\] For} this {we} consider two cases:   {(I) $(x_i)$ does not end with $1^\f$; (II) $(x_i)$ ends with $1^\f$}.

Case (I).  {$(x_i)=\Phi_t(\la_*)$} does not end with $1^\f$. {Then there exists a subsequence $(n_j)\subset {\mathbb{N}}$ such that $x_{n_j}<1$ for all $j\geq 1$.}  By Proposition \ref{prop:monotonicity-Phi-t} there exists a $\de>0$ such that $(\la_*-\de, \la_*+\de)\subset(0,1/3)$, and $\Phi_t$ is monotonic in $(\la_*-\de, \la_*)\cap\La(t)$ and $(\la_*, \la_*+\de)\cap\La(t)$, respectively. Furthermore, if  $\Phi_t$ is strictly decreasing in $(\la_*,\la_*+\de)\cap\La(t)$, then {so is $\Phi_t$} in $(\la_*-\de, \la_*+\de)\cap\La(t)$.   {By deleting the first {finitely many} terms from $(n_j)$} we can assume that for any $j\ge 1$ the equation
\begin{equation}\label{eq:eta-j-1}
t=\Pi(x_1\ldots x_{n_j}^+ 1^\f, \eta_j)
\end{equation}
determines a unique root $\eta_j\in(\la_*-\de,\la_*)$. Accordingly, for $k\ge 1$ we set
\[
\Delta_{k,j}^+(\la_*, t):=\set{\la\in\La(t)\cap(\la_*-\de, \la_*): \Phi_t(\la)=x_1\ldots x_{n_j}^+ d_1d_2\ldots\textrm{ with } d_{ik}\ne {-1} ~\forall i\ge 1}.
\]
Then by Proposition \ref{prop:monotonicity-Phi-t} it follows that
\begin{equation}\label{eq:lower-local-lambda-1}
\Delta_{k,j}^+(\la_*, t)\subset(\eta_j,\la_*)\cap\La(t)~\forall k\ge 1;\quad\textrm{and}\quad \eta_j\nearrow \la_*~\textrm{ {as} }j\to\f.
\end{equation}

If $\Phi_t$ is strictly increasing in $(\la_*,\la_*+\de)\cap\La(t)$, then we can assume that for any $j\ge 1$ {(\ref{eq:eta-j-1})} determines a unique $\tilde\eta_j\in(\la_*,\la_*+\de)$. Accordingly, for $k\ge 1$ we let
\[
\tilde\Delta_{k,j}^+(\la_*, t):=\set{\la\in\La(t)\cap(\la_*, \la_*+\de): \Phi_t(\la)=x_1\ldots x_{n_j}^+ d_1d_2\ldots\textrm{ with } d_{ik}\ne {-1} ~\forall i\ge 1}.
\]
Therefore,
\begin{equation}\label{eq:lower-local-lambda-2}
\tilde \Delta_{k,j}^+(\la_*, t)\subset(\la_*, \tilde\eta_j)\cap\La(t)~\forall k\ge 1;\quad\textrm{and}\quad \tilde\eta_j\searrow \la_*~\textrm{ {as} }j\to\f.
\end{equation}
\begin{lemma}
  \label{lem:localdim-lambda-lower+}
  Let $t\in(0,1)\setminus\set{1/3}$ and $\la_*\in\La(t)\setminus\set{1/3}$. Suppose $\Phi_t(\la_*)$ does not end with $1^\f$.
  \begin{enumerate}[{\rm(i)}]
  \item If $\Phi_t$ is strictly decreasing in $(\la_*,\la_*+\de)\cap\La(t)$, then for  any $k, j\in\N$ there exists   $C_1>0$ such that
  \[
  |\pi_{\eta_j}(\Phi_t(\la_1))-\pi_{\eta_j}(\Phi_t(\la_2))|\le C_1|\la_1-\la_2|
  \]
  for any $\la_1, \la_2\in\Delta_{k,j}^+(\la_*,t)$.

  \item If $\Phi_t$ is strictly increasing in $(\la_*,\la_*+\de)\cap\La(t)$, then for  any $k, j\in\N$ there exists   $C_2>0$ such that
  \[
  |\pi_{\la_*}(\Phi_t(\la_1))-\pi_{\la_*}(\Phi_t(\la_2))|\le C_2|\la_1-\la_2|
  \]
  for any $\la_1, \la_2\in\tilde\Delta_{k,j}^+(\la_*,t)$.
  \end{enumerate}
\end{lemma}
\begin{proof}
Since the proofs of (i) and (ii) are similar, we only prove (i).
Suppose {$\Phi_t$} is strictly decreasing in $(\la_*,\la_*+\de)\cap\La(t)$. Then so is $\Phi_t$ in {$[\eta_j, \lambda_*]\cap \Lambda(t)$} for any $j\geq1$. By a similar argument as in the proof of Lemma \ref{lem:lip-2} one can verify that {for any $\la_1,\la_2\in\Delta_{k,j}^+(\la_*, t)$
\[
|\pi_{\eta_j}(\Phi_t(\la_1))-\pi_{\eta_j}(\Phi_t(\la_2))|\leq \frac{6}{\eta_j^k} |\la_1-\la_2|,
\]
proving (i)}.
\end{proof}

Case (II). {$(x_i)=\Phi_t(\la_*)$} ends with $1^\f$. Then there exists a subsequence $(n_j)\subset\N$ such that $x_{n_j}>-1$ for all $j\ge 1$. The proof is similar to Case (I). By Proposition \ref{prop:monotonicity-Phi-t} there exists a $\de>0$ such that $(\la_*-\de, \la_*+\de)\subset(0,1/3)$, and $\Phi_t$ is monotonic in $(\la_*-\de, \la_*)\cap\La(t)$ and $(\la_*, \la_*+\de)\cap\La(t)$, respectively. Furthermore, if  $\Phi_t$ is strictly decreasing in $(\la_*,\la_*+\de)\cap\La(t)$, then $\Phi_t$ is also strictly decreasing in $(\la_*-\de, \la_*+\de)\cap\La(t)$. So, {by deleting the first {finitely many} terms {from} $(n_j)$} we can assume that for any $j\ge 1$ the equation
\begin{equation}\label{eq:eta-j}
t=\Pi(x_1\ldots x_{n_j}^- (-1)^\f, \ga_j)
\end{equation}
determines a unique root $\ga_j\in(\la_*,\la_*+\de)$. Accordingly, for $k\ge 1$ {let}
\[
\Delta_{k,j}^-(\la_*, t):=\set{\la\in\La(t)\cap(\la_*, \la_*+\de): \Phi_t(\la)=x_1\ldots x_{n_j}^- d_1d_2\ldots\textrm{ with } d_{ik}\ne {-1} ~\forall i\ge 1}.
\]
Then by Proposition \ref{prop:monotonicity-Phi-t} it follows that
\begin{equation}\label{eq:lower-local-lambda-3}
{\Delta_{k,j}^-(\la_*, t)}\subset(\la_*,\ga_j)\cap\La(t)~\forall k\ge 1;\quad\textrm{and}\quad \ga_j\searrow \la_*~\textrm{ {as} }j\to\f.
\end{equation}

If $\Phi_t$ is strictly increasing in $(\la_*,\la_*+\de)\cap\La(t)$, then we can assume that for any $j\ge 1$ {(\ref{eq:eta-j})} determines a unique $\tilde\ga_j\in(\la_*-\de,\la_*)$. Accordingly, for $k\ge 1$ we {set}
\[
\tilde\Delta_{k,j}^-(\la_*, t):=\set{\la\in\La(t)\cap(\la_*-\de, \la_*): \Phi_t(\la)=x_1\ldots x_{n_j}^- d_1d_2\ldots\textrm{ with } d_{ik}\ne {-1} ~\forall i\ge 1}.
\]
Therefore,
\begin{equation}\label{eq:lower-local-lambda-4}
{
\tilde \Delta_{k,j}^-(\la_*, t)\subset(\tilde\ga_j, \la_*)\cap\La(t)~\forall k\ge 1;\quad\textrm{and}\quad \tilde\ga_j\nearrow \la_*~\textrm{ {as} }j\to\f.}
\end{equation}

{By {a} similar augment as in the proof of Lemma \ref{lem:lip-2} we have the following {Lipschitz} property.}
\begin{lemma}
  \label{lem:localdim-lambda-lower-}
  Let $t\in(0,1)\setminus\set{1/3}$ and $\la_*\in\La(t)\setminus\set{1/3}$. Suppose $\Phi_t(\la_*)$   ends with $1^\f$.
  \begin{enumerate}[{\rm(i)}]
  \item If $\Phi_t$ is strictly decreasing in $(\la_*,\la_*+\de)\cap\La(t)$, then for  any $k, j\in\N$ there exists   $C_1>0$ such that
  \[
  |\pi_{\la_*}(\Phi_t(\la_1))-\pi_{\la_*}(\Phi_t(\la_2))|\le C_1|\la_1-\la_2|
  \]
  for any $\la_1, \la_2\in\Delta_{k,j}^-(\la_*,t)$.

  \item If $\Phi_t$ is strictly increasing in $(\la_*,\la_*+\de)\cap\La(t)$, then for  any $k, j\in\N$ there exists   $C_2>0$ such that
  \[
  |\pi_{\tilde\ga_j}(\Phi_t(\la_1))-\pi_{\tilde\ga_j}(\Phi_t(\la_2))|\le C_2|\la_1-\la_2|
  \]
  for any $\la_1, \la_2\in\tilde\Delta_{k,j}^-(\la_*,t)$.
  \end{enumerate}
\end{lemma}

\begin{proof}[Proof of Proposition \ref{th:local-dim-lambda} ({a lower bound})]
Let $t\in(0,1)\setminus\set{1/3}$ and $\la_*\in\La(t)\setminus\set{1/3}$. If $\Phi_t(\la_*)$ does not end with $1^\f$ and $\Phi_t$ is strictly decreasing in $(\la_*,\la_*+\de)\cap\La(t)$ for some $\de>0$, then by (\ref{eq:lower-local-lambda-1})  and Lemma \ref{lem:localdim-lambda-lower+} (i) it follows that for any $j\geq1$,
\begin{align*}
\dim_H(\La(t)\cap(\la_*-\de, \la_*+\de))&\ge \dim_H\pi_{\eta_j}(\Phi_t(\Delta_{k,j}^+(\la_*,t)))\\
&\ge \frac{(k-1)\log 3+\log 2}{-k\log\eta_j} \to \frac{\log3}{-\log\eta_j} {\quad \textrm{as } k\to\infty}.
\end{align*}
Since $\eta_j\nearrow \la_*$ as $j\to\f$ by (\ref{eq:lower-local-lambda-1}),  we conclude  that
\begin{equation}\label{eq:local-lower-bound}
\lim_{\de\to 0^+}\dim_H(\La(t)\cap(\la_*-\de, \la_*+\de))\ge \lim_{j\to\f} \frac{\log 3}{{-\log \eta_j}}=\frac{\log 3}{-\log\la_*}.
\end{equation}

Similarly, if $\Phi_t(\la_*)$ does not end with $1^\f$ and $\Phi_t$ is strictly increasing in $(\la_*, \la_*+\de)\cap\La(t)$, then by (\ref{eq:lower-local-lambda-2}), Lemma \ref{lem:localdim-lambda-lower+} (ii) and the same argument as above we can prove (\ref{eq:local-lower-bound}). Moreover, if $\Phi_t(\la_*)$ ends with $1^\f$, then (\ref{eq:local-lower-bound}) can be proved by a similar argument as above together with (\ref{eq:lower-local-lambda-3}), (\ref{eq:lower-local-lambda-4}) and Lemma \ref{lem:localdim-lambda-lower-}.
\end{proof}

 \subsection{{{An upper bound {on} the local dimension of $\Lambda(t)$}}}  Now we turn to prove that
 {\[\lim_{\de\to 0^+}\dim_H(\La(t)\cap(\la_*-\de, \la_*+\de))\le\frac{\log 3}{-\log\la_*}\quad\forall \la_*\in\La(t)\setminus\set{1/3}.\]}In view of our proof of Proposition \ref{prop:monotonicity-Phi-t} we consider $t\in(0,1/9]$, {$t\in(1/9,4/27)$ and} $t\in[4/27,1/3)\cup(1/3,1)$, separately.

\begin{lemma}
  \label{lem:upper-local-lambda-4/27-1}
 {Let} $t\in[4/27, 1/3)\cup(1/3,1)${. Then} for any $\theta\in(\min\La(t), 1/3)$  there exists $C>0$ such that
  \[
  |\pi_{\theta}(\Phi_t(\la_1))-\pi_{\theta}(\Phi_t(\la_2))|\ge C|\la_2-\la_1|
  \]
  for any $\la_1, \la_2\in\La(t)\cap(0, \theta]$.
\end{lemma}
\begin{proof}
First {we consider} $t\in[4/27, 1/3)$. {Take} $\theta \in(\min\La(t),1/3)$. {By} the proof of Lemma \ref{lem:mon-Phi-t-4/27-1} it follows that $\Phi_t(\la)\in[01(-1)0^\f, 01^\f]$ for any $\la\in\La(t)\cap(0,\theta]$. {Furthermore, by} the proof of Lemma \ref{lem:monotonicity-lambda} (iii) {we obtain} that
\begin{equation}
  \label{eq:upper-local-21}
  \begin{split}
  \Pi_2((i_n),\la)&\ge\min\set{1-4\la+3\la^2, 1-2\la-3\la^2} \ge 1-4\theta+3\theta^2=:C_1>0
  \end{split}
\end{equation}
for any $(i_n)\in [01(-1)0^\f, 01^\f]$ and $\la\in\La(t)\cap(0,\theta]$.

  Now take $\la_1, \la_2\in\La(t)\cap(0, \theta]$ with $\la_1<\la_2$. {Then} $(i_n):=\Phi_t(\la_1)\succ\Phi_t(\la_2)=:(j_n)$ by  Lemma \ref{lem:mon-Phi-t-4/27-1}. {So} there exists $N\in\N_{\ge 3}$ such that $i_1\ldots i_{N-1}=j_1\ldots j_{N-1}$ and $i_N>j_N$.   {Thus,}
  \begin{equation}
    \label{eq:upper-local-23}
     \begin{split}
      |\pi_{\theta}(\Phi_t(\la_1))-\pi_{\theta}(\Phi_t(\la_2))| &=(1-\theta)\sum_{n=N}^\f(i_n-j_n)\theta^{n-1}\\
      & \ge (1-\theta)\theta^{N-1}-(1-\theta)\sum_{n=N+1}^\f 2\theta^{n-1} =(1-3\theta)\theta^{N-1}.
      \end{split}
  \end{equation}

  {On} the other hand, note that
  $
  (1-\la_1)\sum_{n=1}^\f i_n\la_1^{n-1}=t=(1-\la_2)\sum_{n=1}^\f j_n\la_2^{n-1}.
 $
  Then by (\ref{eq:upper-local-21}) it follows that
  \begin{align*}
    &\quad(1-\la_1)\sum_{n=N}^\f(i_n-1)\la_1^{n-1}-(1-\la_2)\sum_{n=N}^\f (j_n-1)\la_2^{n-1}\\
    &=\left((1-\la_2)\sum_{n=1}^{N-1} j_n\la_2^{n-1}+(1-\la_2)\sum_{n=N}^\f\la_2^{n-1}\right)\\
    &\qquad -\left((1-\la_1)\sum_{n=1}^{N-1}j_n\la_1^{n-1}+(1-\la_1)\sum_{n=N}^\f\la_1^{n-1}\right)\\
    &=\Pi(j_1\ldots j_{N-1}1^\f, \la_2)-\Pi(j_1\ldots j_{N-1}1^\f, \la_1)\\
    &=\Pi_2(j_1\ldots j_{N-1}1^\f, \la_*)(\la_2-\la_1)\ge C_1(\la_2-\la_1)
  \end{align*}
  for some $\la_*\in[\la_1,\la_2]$, {where the third equality follows by the mean value theorem. So,}
  \begin{equation}
    \label{eq:upper-local-24}
    \begin{split}
    \la_2-\la_1&\le \frac{1}{C_1}\left[(1-\la_1)\sum_{n=N}^\f 2\la_1^{n-1}+(1-\la_2)\sum_{n=N}^\f 2\la_2^{n-1}\right]\\
    &\le\frac{2}{C_1}(\la_1^{N-1}+\la_2^{N-1})\le\frac{4}{C_1}\theta^{N-1}.
    \end{split}
  \end{equation}
  Hence, by (\ref{eq:upper-local-23}) and (\ref{eq:upper-local-24}) we conclude that
  \[
  |\pi_\theta(\Phi_t(\la_1))-\pi_\theta(\Phi_t(\la_2))|\ge \frac{C_1(1-3\theta)}{4}{|\la_1-\la_2|}
  \]
 as required.

 Next we consider $t\in(1/3,1)$. Then by the proof of Lemma \ref{lem:mon-Phi-t-4/27-1} there exists $k\in\N$ such that $\Phi_t(\la)\in[1(-1)^\f, 1^k(-1)^\f]$ for any $\la\in\La(t)\cap(0,\theta]$. Since $\min\La(t)=(1-t)/2$, by (\ref{eq:derivative}) it follows that
 \[
 |\Pi_2((i_n),\la)|\ge|\Pi_2(1^k(-1)^\f,\la)|= 2k\la^{k-1}\ge2k\left(\frac{1-t}{2}\right)^{k-1}=:C_2>0.
 \]
 By a similar argument as above one can prove   that for any $\la_1, \la_2\in\La(t)\cap(0,\theta]$
 \[
 |\pi_\theta(\Phi_t(\la_1))-\pi_\theta(\Phi_t(\la_2))|\ge \frac{C_2(1-3\theta)}{4}|\la_1-\la_2|,
 \]
 completing the proof.
\end{proof}

Recall from Lemma \ref{lem:mon-Phi-t-1} that  {$\la_\diamond=\lambda_\diamond(t) \in(0,1/3)$} satisfies  $\la_\diamond(1-\la_{\diamond})^2=t$. Furthermore, by (\ref{eq:la-diamond-tau}) we have $t<\la_\diamond(t)<\sqrt{t}\le 1/3$ for all $t\in(0,1/9]$.
\begin{lemma}
   \label{lem:upper-local-lambda-0-1/9}
  Let $t\in(0,1/9]$.
  \begin{enumerate}[{\rm(i)}]
  \item For any $\theta\in[t, \la_\diamond]$  there exists $C>0$ such that for any $\la_1, \la_2\in\La(t)\cap(0, \la_\diamond]$
  \[
  |\pi_{\theta}(\Phi_t(\la_1))-\pi_{\theta}(\Phi_t(\la_2))|\ge C|\la_2-\la_1|;
  \]

  \item  For any $\theta\in(\la_\diamond, \sqrt{t})$  there exists $C'>0$ such that for any $\la_1,\la_2\in\La(t)\cap(\la_\diamond, \sqrt{t})$
  \[
  |\pi_{\theta}(\Phi_t(\la_1))-\pi_{\theta}(\Phi_t(\la_2))|\ge C'|\la_1-\la_2|;
  \]

  \item For any {$\theta\in[\sqrt{t}, 1/3)$}  there exists $C''>0$ such that for any $\la_1, \la_2\in\La(t)\cap[\sqrt{t}, \theta]$
  \[
  |\pi_{\theta}(\Phi_t(\la_1))-\pi_{\theta}(\Phi_t(\la_2))|\ge C''|\la_2-\la_1|.
  \]
  \end{enumerate}
\end{lemma}
\begin{proof}
For (i), note  {by the proof of Lemma \ref{lem:mon-Phi-t-1}} that $\Phi_t(\la)\in[01(-1)0^\f, 01^\f]$ for any $\la\in\La(t)\cap(0,\la_\diamond]$. Then by (\ref{eq:upper-local-21}) it follows that
$
\Pi_2((i_n), \la)\ge 1-4\la_\diamond+3\la_{\diamond}^2>0
$
for all $(i_n)\in[01(-1)0^\f, 01^\f]$ and $\la\in(0,\la_\diamond]$. So, (i)  follows by the same argument as in the proof of Lemma
 \ref{lem:upper-local-lambda-4/27-1}.

For (ii), {note by (\ref{eq:dec14-2})}  that {for any $\la\in\La(t)\cap(\la_\diamond,\sqrt{t})$ we have
$\Phi_t(\la)\in(001^\f, 01(-1)0^\f)$, and the proof of Lemma \ref{lem:mon-phi-t-0-1/9} yields that $\la<\la_{\Phi_t(\la)}$. Then $\Pi_2(\Phi_t(\la),\la)>0$ for all $\la\in\La(t)\cap(\la_\diamond, \sqrt{t})$. Furthermore,} observe by Lemma \ref{lem:monotonicity-lambda} that $\Pi_2(\Phi_t(\la_{\diamond}), \la_{\diamond})=\Pi_2(01(-1)0^\f, \la_{\diamond})>0$ and $\Pi_2(\Phi_t(\sqrt{t}), \sqrt{t})=\Pi_2(001^\f, \sqrt{t})>0$. Since the map $t\mapsto\Pi_2(\Phi_t(\la),\la)$ is continuous by Lemma \ref{lem:continuity-Phi-t} {and $\La(t)$ is compact by Proposition \ref{prop:topology-La(t)}}, it follows that
\[
C_1:=\inf_{\la\in[\la_\diamond, \sqrt{t}]\cap\La(t)}\Pi_2(\Phi_t(\la),\la)=\min_{\la\in[\la_\diamond, \sqrt{t}]\cap\La(t)}\Pi_2(\Phi_t(\la),\la)>0.
\]
Now take $\la_1,\la_2\in\La(t)\cap(\la_\diamond, \sqrt{t})$ with $\la_1<\la_2$. Then by Lemma \ref{lem:mon-phi-t-0-1/9} we have $(i_n):=\Phi_t(\la_1)\succ\Phi_t(\la_2)=(j_n)$. So, there exists $N\in\N_{\ge 4}$ such that $i_1\ldots i_{N-1}=j_1\ldots j_{N-1}$ and $i_N>j_N$. By the same argument as in the proof of Lemma \ref{lem:upper-local-lambda-4/27-1} one can show that
\[
|\pi_{\theta}(\Phi_t(\la_1))-\pi_{\theta}(\Phi_t(\la_2))|\ge (1-3\theta)\theta^{N-1},
\]
and there exists $\la_*\in[\la_1,\la_2]$ such that
\begin{align*}
  &(1-\la_1)\sum_{n=N}^\f (i_n-1)\la_1^{n-1}-(1-\la_2)\sum_{n=N}^\f(j_n-1)\la_2^{n-1}\\
  =&\Pi(j_1\ldots j_{N-1}1^\f,\la_2)-\Pi(j_1\ldots j_{N-1}1^\f,\la_1)\\
  =&\Pi_2(j_1\ldots j_{N-1}1^\f, \la_*)(\la_2-\la_1)\\
  \ge&\Pi_2(j_1\ldots j_{N-1}1^\f, \la_2)(\la_2-\la_1)\\
  \ge &\Pi_2(\Phi_t(\la_2), \la_2)(\la_2-\la_1)\ge C_1(\la_2-\la_1),
\end{align*}
where the first inequality follows by the concavity of {$\Pi(j_1\ldots j_{N-1}1^\infty,\cdot)$} and the second inequality follows by Lemma \ref{lem:monotonicity-in}.
Therefore, one can deduce from this that
\[
 |\pi_{\theta}(\Phi_t(\la_1))-\pi_{\theta}(\Phi_t(\la_2))|\ge \frac{C_1(1-3\theta)}{4}|\la_1-\la_2|,
\]
establishing (ii).

For (iii), by the proof of {Lemma \ref{lem:mon-phi-t-0-1/9}} it follows that
$
\Phi_t(\la)\in[\Phi_t(1/3), \Phi_t(\sqrt{t})]\subset(0^\f, 001^\f]$  for all $\la\in\La(t)\cap[\sqrt{t}, 1/3].
$
So, by the proof of Lemma \ref{lem:monotonicity-lambda} (i) there exists $k\in\N_{\ge 3}$ such that
\begin{align*}
\Pi_2((i_n),\la)&\ge \Pi_2(\Phi_t(1/3),\la)\ge \la^{k-2}(k-1-2k\la)\ge \min_{\la\in[\sqrt{t},\theta]}\la^{k-1}(k-1-2k\la){=:C_2>0}.
\end{align*}
for all $(i_n)\in[\Phi_t(1/3), \Phi_t(\sqrt{t})]$ {and $\lambda\in[\sqrt{t},\theta]$.} Thus {(iii)} follows by
the same argument as in the proof of Lemma \ref{lem:upper-local-lambda-4/27-1}.
\end{proof}

Note by the proof of Proposition \ref{prop:monotonicity-Phi-t} (ii) that {$\tau=\tau(t) \in[1/4, 1/3)$ for $t\in(1/9,4/27)$. Furthermore, by (\ref{eq:la-diamond-tau}) we have}  $t<\la_\diamond(t)<\tau(t)<1/3$ for all {$t\in(1/9,4/27)$}.

\begin{lemma}
   \label{lem:upper-local-lambda-1/9-4/27}
  Let $t\in(1/9, 4/27)$.
  \begin{enumerate}[{\rm(i)}]
   \item For any $\theta\in[t, \la_\diamond]$  there exists $C>0$ such that for any $\la_1, \la_2\in\La(t)\cap(0, \la_\diamond]$
  \[
  |\pi_{\theta}(\Phi_t(\la_1))-\pi_{\theta}(\Phi_t(\la_2))|\ge C|\la_2-\la_1|;
  \]
  \item For any $\theta\in(\la_\diamond, \tau)$  there exists $C'>0$ such that for any $\la_1, \la_2\in\La(t)\cap[\la_\diamond, \theta]$
  \[
  |\pi_{\theta}(\Phi_t(\la_1))-\pi_{\theta}(\Phi_t(\la_2))|\ge C'|\la_2-\la_1|;
  \]

  \item    For any $\theta\in(\tau, 1/3)$  there exists $C''>0$ such that for any $\la_1, \la_2\in\La(t)\cap[\theta, 1/3]$
  \[
  |\pi_{\theta}(\Phi_t(\la_1))-\pi_{\theta}(\Phi_t(\la_2))|\ge C''|\la_2-\la_1|.
  \]

  \end{enumerate}
\end{lemma}
\begin{proof}
The proof is similar to that for Lemma \ref{lem:upper-local-lambda-0-1/9}. More precisely,
(i) follows by the same argument as in the proof of Lemma \ref{lem:upper-local-lambda-0-1/9} (i). In view of the proofs of Lemmas \ref{lem:mon-phi-t-1/9-1/8} and \ref{lem:mon-phi-t-1/8-4/27}, (ii) and (iii) follow by a similar argument as in the proof of Lemma \ref{lem:upper-local-lambda-0-1/9} (ii).
\end{proof}

\begin{proof}
  [Proof of Proposition \ref{th:local-dim-lambda} ({an upper bound})]
  Let $\la_*\in\La(t)\setminus\set{1/3}$. We   consider  the following two cases.

  Case I. $t\in[4/27, 1/3)\cup(1/3,1)$.  Then by Lemma \ref{lem:upper-local-lambda-4/27-1} it follows that for any $\de\in(0,1/3-\la_*)$,
  \[
  \dim_H(\La(t)\cap(\la_*-\de, \la_*+\de))\le \dim_H\pi_{\la_*+\de}(\Phi_t(\La(t)\cap(\la_*-\de, \la_*+\de)))\le \frac{\log 3}{-\log(\la_*+\de)},
  \]
{which implies}
  \begin{equation}
    \label{eq:upper-local-25}
    \lim_{\de\to 0^+}\dim_H(\La(t)\cap(\la_*-\de, \la_*+\de))\le\frac{\log 3}{-\log\la_*}
  \end{equation}
  as desired.

  Case II. $t\in(0,4/27)$.
  Then by Lemmas \ref{lem:upper-local-lambda-0-1/9} and \ref{lem:upper-local-lambda-1/9-4/27} it follows that for any $\de\in(0,1/3-\la_*)$
  \[
  \dim_H(\La(t)\cap(\la_*-\de, \la_*))\le\frac{\log 3}{-\log\la_*}\quad\textrm{and}\quad \dim_H(\La(t)\cap(\la_*,\la_*+\de))\le\frac{\log 3}{-\log(\la_*+\de)}.
  \]
  Letting $\de\to 0^+$ we also obtain (\ref{eq:upper-local-25}).
\end{proof}

\begin{proof}[Proof of Theorem \ref{thm: La(t)}]
By Proposition  \ref{prop:topology-La(t)}, Corollary \ref{cor:full-dimension} and Proposition \ref{th:local-dim-lambda} it suffices to prove that $\La(t)$ has zero Lebesgue measure for all $t\in(0,1)\setminus\set{1/3}$. Note by Proposition \ref{th:local-dim-lambda} that for each $\la\in\La(t)\setminus\set{1/3}$ there exists $\de_\la>0$ such that $\dim_H(\La(t)\cap(\la-\de_\la,\la+\de_\la))<1$, and thus $\La(t)\cap(\la-\de_\la,\la+\de_\la)$ has zero Lebesgue measure. Since by Proposition \ref{prop:topology-La(t)} that $\La(t)$ is compact, for any $n\in\N_{\ge 4}$ the segment $\La(t)\cap(0,1/3-1/n]$ can be covered by {finitely many open intervals}, say $\set{(\la_i-\de_{\la_i}, \la_i+\de_{\la_i}): i=1,\ldots, N}$. {This implies} that $\La(t)\cap(0,1/3-1/n]$ has zero Lebesgue measure for all $n\ge 4$. Therefore, $\La(t)=\set{1/3}\cup\bigcup_{n=4}^\f(\La(t)\cap(0,1/3-1/n])$ {is a Lebesgue null set}.
\end{proof}

\section{Intersections with different Hausdorff and packing dimensions}\label{sec:level-set-not}

In this section {we consider} the set
\[
\Lambda_{not}(t):=\{\lambda\in\Lambda(t):\dim_H(C_\lambda\cap(C_\lambda+t))\ne \dim_P(C_\lambda\cap(C_\lambda+t))\},
\]
and show that it has full Hausdorff dimension. {Our proof is motivated by the work of \cite{Sergio_Mykola_Drygoriy_2005} to construct} subsets of $\Lambda_{not}(t)$ {with} Hausdorff {dimension} arbitrarily  close  to $1$. {Recall by (\ref{eq:dim-intersection}) that
\[
\dim_H(C_\lambda\cap(C_\lambda+t))=\frac{\log2}{-\log\lambda}\underline {freq}_0(\Phi_t(\la)),\quad \dim_P(C_\la\cap(C_\la+t))=\frac{\log 2}{-\log\la}\overline {freq}_0(\Phi_t(\la)),
\]
where for $(i_n)\in\set{-1,0,1}^\N$ we {denote}
\[\underline {freq}_0((i_n)):=\liminf_{n\to\infty}\frac{\#\{1\leq k\leq n:i_k=0\}}{n},\quad \overline {freq}_0((i_n)):=\limsup_{n\to\infty}\frac{\#\{1\leq k\leq n:i_k=0\}}{n}.\]
Then $\La_{not}(t)$ can be rewritten as
\begin{equation}\label{eq:Lambda-not}
{\Lambda_{not}(t)=\set{\lambda\in\Lambda(t):\underline {freq}_0(\Phi_t(\la))\ne \overline {freq}_0(\Phi_t(\la))}.}
\end{equation}}
\begin{proposition}\label{prop:Lambda-not}
For any $t\in(0,1)\setminus \{1/3\}$ we have {$\dim_H \Lambda_{not}(t)=1$.}
\end{proposition}

Similar to the proof of Proposition \ref{prop:local-dim-1/3}, we consider two cases: $t\in(0,1/9]\cup[4/27,1/3)$ and $t\in(1/9,4/27)\cup(1/3,1)$. {Note by Proposition \ref{prop:monotonicity-Phi-t} that the map $\Phi_t$ is strictly monotonic for $t\in(0,1/9]\cup[4/27,1/3)$ but piecewise monotonic for $t\in(1/9,4/27)\cup(1/3,1)$.} Since the proofs of the two cases are similar,  in the following we only
consider $t\in(1/9,4/27)\cup(1/3,1)$. By Proposition \ref{prop:monotonicity-Phi-t} it follows that there exists a small $\de>0$ such that   $\Phi_t$ is strictly increasing in $\Lambda(t)\cap[1/3-\de, 1/3]$.
Let $(x_i)=\Phi_t(1/3)$. Then $(x_i)$  does not end with $(-1)^\infty$. So there exists a {subsequence} $\{n_j\}\subset\mathbb{N}$  such that $x_{n_j}>-1$  for all $j\ge 1$.
By deleting the first {finitely many} terms {from} $(n_j)$  we may assume that for each $j\ge 1$ the equation
\[
t=\Pi(x_1\ldots x_{n_j}^-(-1)^\f, \ga_j)
\]
determines a unique $\gamma_j\in(1/3-\de,1/3)$.
{Clearly,} $\ga_j\nearrow 1/3$ as $j\to \f$.

Now, for $q\in\N$ let $\Sigma_{q,j}(t)$ consists of all $\la\in\La(t)\cap(1/3-\de, 1/3)$ such that
\begin{equation}\label{eq:Sigma-qj}
\begin{split}
  \Phi_t(\la)=~&x_1\ldots x_{n_j}^-d_1d_2\ldots\\
  =~&x_1\ldots x_{n_j}^-\\
  &d_{r_0+1}d_{r_0+2}\ldots d_{r_0+3q}110\\
  &d_{r_1+1}d_{r_1+2}\ldots d_{r_1+2\cdot 3q}111100\\
  &\cdots\\
  &d_{r_{m-1}+1}d_{r_{m-1}+2}\ldots d_{r_{m-1}+2^{m-1}\cdot 3q}1^{2^{m}}0^{2^{m-1}}\\
  &d_{r_m+1}d_{r_m+2}\ldots d_{r_m+2^m\cdot 3q}1^{2^{m+1}}0^{2^m}\\
  &\cdots,
\end{split}
\end{equation}
where
\[r_m:=3(q+1)(1+2+\cdots+2^{m-1})=3(q+1)(2^m-1)\quad\textrm{and}\quad d_{r_m+kq}\ne -1\]
 for all $k\in\set{1,\ldots, {2^m\cdot 3}}$ and $m\ge 0$. Here we point out that $r_0=0$.
Then by Proposition \ref{prop:monotonicity-Phi-t}   it follows that
$
\Sigma_{q,j}(t)\subset\La(t)\cap(\ga_j,1/3)$ for all $q\ge 1.
$
In the following we show that $\Sigma_{q, j}(t)$ is {a} subset of $\La_{not}(t)$.

\begin{lemma}\label{lem:Lambda-not-1}
Let $t\in(1/9,4/27)\cup(1/3,1)$. Then
 \[\Sigma_{q,j}(t)\subset\La_{not}(t)\cap(\ga_j, 1/3)\quad\forall q, j\in\N.\]
\end{lemma}

\begin{proof}
Take $\lambda\in\Sigma_{q,j}(t)$. Then $\Phi_t(\la)=x_1\ldots x_{n_j}^-d_1d_2\ldots$ {is defined as} in (\ref{eq:Sigma-qj}). By (\ref{eq:Lambda-not}) it suffices to prove that the frequency of digit zero in the sequence $(d_i)$ does not exist. For $n\in\N$ let $N_0((d_i), n)$ be the number of   zeros in the word $d_1\ldots d_n$. Then
\[{
N_0((d_i), r_m)=\xi_m((d_i))+(1+2+\cdots+2^{m-1})=\xi_m((d_i))+2^m-1,}
\]
where $\xi_m((d_i))$ denotes the number of zeros in the word
\[{
d_{r_0+1}d_{r_0+2}\ldots d_{r_0+3q}\;d_{r_1+1}d_{r_1+2}\ldots d_{r_1+2\cdot 3q}\;\cdots \; d_{r_{m-1}+1}d_{r_{m-1}+2}\ldots d_{r_{m-1}+2^{m-1}\cdot 3q}.}\]
This implies that
\begin{equation}
  \label{eq:limit-1}
  \lim_{m\to\f}\frac{N_0((d_i), r_m)}{r_m}=\lim_{m\to\f}\frac{\xi_m((d_i))+2^m-1}{3(q+1)(2^m-1)}=\frac{1}{3(q+1)}\left(1+ \lim_{m\to\f}\frac{\xi_m((d_i))}{2^m}\right).
\end{equation}

Similarly, for $\ell_m:=r_m-2^{m-1}$ we have
\[
N_0((d_i), \ell_m)=\xi_m((d_i))+(1+2+\cdots+2^{m-2})=\xi_m((d_i))+2^{m-1}-1,
\]
and thus
\begin{equation}
  \label{eq:limit-2}
  \begin{split}
  \lim_{m\to\f}\frac{N_0((d_i), \ell_m)}{\ell_m}&=\lim_{m\to\f}\frac{\xi_m((d_i))+2^{m-1}-1}{3(q+1)(2^m-1)-2^{m-1}}\\
  &=\frac{1}{6q+5}\left(1+2\lim_{m\to\f}\frac{\xi((d_i),m)}{2^m}\right).
  \end{split}
\end{equation}
Combining    (\ref{eq:limit-1}) with (\ref{eq:limit-2}), if the limit $\lim_{m\to\f}\frac{\xi_m((d_i))}{2^m}$  does not exist, then both limits  $\lim_{m\to\f}\frac{N_0((d_i), r_m)}{r_m}$ and $\lim_{m\to\f}\frac{N_0((d_i), \ell_m)}{\ell_m}$ do not exist, and thus the frequency of digit zero in $(d_i)$ does not exist. If the limit $\lim_{m\to\f}\frac{\xi_m((d_i))}{2^m}$ {exists}, call this limit $\zeta$, then since $\xi_m((d_i))\le 3q(1+2+\ldots+2^{m-1})=3q(2^m-1)$, we have $\zeta\le 3q$. Therefore, by  (\ref{eq:limit-1}) and (\ref{eq:limit-2}) it follows that
\[\lim_{m\to\f}\frac{N_0((d_i), r_m)}{r_m}=\frac{\zeta+1}{3(q+1)}>\frac{2\zeta+1}{6q+5}=\lim_{m\to\f}\frac{N_0((d_i), \ell_m)}{\ell_m},\]
 which again implies that the frequency of digit zero in $(d_i)$ does not exist. So, $\la\in\La_{not}(t)$, completing  the proof.
\end{proof}

  Next we  give {a lower bound} for the Huasdorff dimension of $\Sigma_{q,j}(t)$.

\begin{lemma}\label{lem:Lambda-not-2}
Let $t\in(1/9,4/27)\cup(1/3,1)$. Then for any $q, j\in\N$ we have
\[dim_H\Sigma_{q,j}(t) \geq  \frac{ (q-1)\log3+ \log2}{- (q+1)\log\gamma_j}.\]
\end{lemma}
\begin{proof}
  Note by Lemma \ref{lem:Lambda-not-1} that $\Sigma_{q,j}(t)\subset\Lambda_{not}(t)\cap(\ga_j,1/3)$. Then for any $\la\in\Sigma_{q,j}(t)$ the sequence $\Phi_t(\la)$ begins with $x_1\ldots x_{n_j}^-$ and the tail sequence does not contain $q$ consecutive $(-1)$s. Then by the same argument as in the proof of Lemma \ref{lem:lip-1} there exists $C>0$ such that
  \[
  |\pi_{\ga_j}(\Phi_t(\la_1))-\pi_{\ga_j}(\Phi_t(\la_2))|\le C|\la_1-\la_2|
  \]
  for any $\la_1,\la_2\in\Sigma_{q,j}(t)$. This implies that
  $\dim_H\Sigma_{q,j}(t)\ge \dim_H\pi_{\ga_j}(\Phi_t(\Sigma_{q,j}(t)))$. So it suffices to prove
  \begin{equation}\label{eq:Lambda-not-1}
  \dim_H\pi_{\ga_j}(\Phi_t(\Sigma_{q,j}(t)))\ge\frac{ (q-1)\log 3+ \log 2}{- (q+1)\log\ga_j}.
  \end{equation}

  Observe by (\ref{eq:Sigma-qj}) that
  \begin{equation}\label{eq:Lambda-not-2}
  \dim_H\pi_{\ga_j}(\Phi_t(\Sigma_{q,j}(t)))=\dim_H\pi_{\ga_j}{\left(\prod_{m=0}^\infty E_m(q)\right)},
  \end{equation}
  where
  \[{
  E_m(q)}:=\left( \set{-1,0,1}^{q-1}\times\set{0,1}\right)^{3\cdot 2^m}\times\set{1^{2^{m+1}}0^{2^m}}.
  \]
  Note that each word in $E_m(q)$ has length {$3(q+1)2^m$}, and $\#E_m(q)=(3^{q-1}\cdot 2 )^{3\cdot 2^m}$. Furthermore, $\prod_{m=0}^\infty E_m(q)$ is the set of infinite sequences by concatenating  words from {each} $E_m(q)$.
  So $\pi_{\ga_j}(\prod_{m=0}^\infty E_m(q))$ is a homogeneous Moran set satisfying the strong separation condition. {Hence,} by \cite[Theorem 2.1]{Feng-Wen-Wu-1997} it follows that
  \begin{align*}
    \dim_H\pi_{\ga_j}{\left({\prod_{m=0}^\infty E_m(q)}\right)}&\ge \liminf_{m\to\f}\frac{\log\prod_{\ell=0}^{m-1}(3^{q-1}\cdot 2)^{3\cdot 2^\ell}}{-{\sum_{l=0}^{m-1}3(q+1)2^l}\log \ga_j }\\
    &=\liminf_{m\to\f}\frac{\sum_{\ell=0}^{m-1}3\cdot 2^\ell\log (3^{q-1}\cdot 2)}{-3(q+1)(2^m-1)\log \ga_j}\\
    &=\liminf_{m\to\f}\frac{3(2^m-1)\log (3^{q-1}\cdot 2)}{-3(q+1)(2^m-1)\log \ga_j}=\frac{ (q-1)\log 3+ \log 2}{- (q+1)\log\ga_j}.
  \end{align*}
  This, together with (\ref{eq:Lambda-not-2}), proves (\ref{eq:Lambda-not-1}).
\end{proof}

\begin{proof}
  [Proof of Proposition \ref{prop:Lambda-not}]
Take $t\in(0,1)\setminus\set{1/3}$. Since the proof for   $t\in(0,1/9]\cup[4/27,1/3)$ is analogous,  we only consider for $t\in(1/9,4/27)\cup(1/3, 1)$. By Lemmas \ref{lem:Lambda-not-1} and \ref{lem:Lambda-not-2} it follows that for any $q\in{\mathbb{N}}$
\[
\dim_H\Lambda_{not}(t)\geq
\dim_H\Sigma_{q,j}(t) \geq  \frac{ (q-1)\log3+ \log2}{- (q+1)\log\gamma_j}\to \frac{(q-1)\log3+\log2}{(q+1)\log3} \quad \textrm{as } j\to\infty.
\]
Letting $q\to\f$ we obtain $\dim_H\Lambda_{not}(t)\ge  1$.
Since the reverse inequality is obvious, this proves $\dim_H\Lambda_{not}(t)=1$.
\end{proof}

\section{Proof of Theorem \ref{th:La(t)-alpha}}\label{sec:Proof-theorem2}
 In this section {we} consider the level set
\[
\Lambda_\beta(t):=\set{\lambda\in\La(t): \dim_H \big(C_\lambda\cap(C_\lambda+t)\big)=\dim_P \big(C_\lambda\cap(C_\lambda+t)\big)=\beta\frac{\log 2}{- \log\la}    }
\]
for $\beta\in[0,1]$, and prove Theorem \ref{th:La(t)-alpha}. {Recall from (\ref{eq:entropy}) the entropy function $h(p_1,p_2,p_3)=-\sum_{i=1}^3 p_i\log p_i$ for a probability vector $(p_1, p_2, p_3)$.}

\begin{proposition}
  \label{th:level-set}
  Let $t\in(0,1)\setminus\set{1/3}$. Then for any $\beta\in[0,1]$
  \[
\dim_H\La_\beta(t)=\dim_P\La_\beta(t)=\frac{h(\frac{1-\beta}{2}, \beta, \frac{1-\beta}{2})}{\log 3}.
\]
\end{proposition}
It is worth   mentioning that the dimension of $\La_\beta(t)$ is independent of $t$.   Since the proof for  $t\in(0,1/9]\cup(4/27, 1/3)$ is  similar,   we only consider $t\in(1/9,4/27)\cup(1/3,1)$.  Then {$(x_i)=\Phi_t(1/3)$} does not end with $(-1)^\f$. So there exists a {subsequence} $(n_j)\subset\N$ such that $x_{n_j}>-1$ for all $j\ge 1$. By Proposition \ref{prop:monotonicity-Phi-t} there exists $\de>0$ such that $\Phi_t$ is strictly increasing in $\La(t)\cap(1/3-\de, 1/3)$. By deleting the first {finitely many} terms from $(n_j)$ we may assume that for each $j\geq 1$ the equation {\[t=\pi_{\ga_j}(x_1\ldots x_{n_j}^-(-1)^\f)\]
determines} a unique $\ga_j\in(1/3-\de, 1/3)$. Then by Lemma \ref{lem:continuity-Phi-t} it follows that  $\ga_j\nearrow 1/3$ as $j\to\f$.

{Given} $\beta\in[0,1]$, we first prove {the lower bound \[\dim_H\Lambda_\beta(t)\geq \frac{h(\frac{1-\beta}{2}, \beta, \frac{1-\beta}{2})}{\log 3}.\]}For $k,j\in{\mathbb{N}}$ let $\F_{k,j}^-(t,\beta)$ consist  of all $\la\in\La(t)\cap(1/3-\de, 1/3)$ such that
{\[
\Phi_t(\la)=x_1\ldots x_{n_j}^-d_1d_2\ldots
\]
satisfying
\begin{equation}\label{eq:freqency-F-kj}
 d_{nk+1}\ldots d_{nk+k} \ne  (-1)^k \;\; {\forall n\ge 0}\quad\textrm{and}\quad
{freq_0((d_i)):=\underline{freq}_0((d_i))=\overline{freq}_0((d_i))=\beta}.
\end{equation}}Then by Proposition \ref{prop:monotonicity-Phi-t} it follows that for each $j\in{\mathbb{N}}$,
\begin{equation}\label{eq:subset-F-kj-alpha-}
\F_{k,j}^-(t,\beta)\subset\La_\beta(t)\cap(\ga_j, 1/3)\quad\forall k\ge 1.
\end{equation}
Note that for any sequence $(c_i)\in \Phi_t(\F_{k,j}^-(t,\beta))$ the tail sequence $c_{n_j+1}c_{n_j+2}\ldots$ does not contain $2k-1$ consecutive $(-1)$s. Then by the same argument  as in the proof of {Lemma \ref{lem:lip-2}} it follows that
 \begin{equation}\label{eq:subset-F-kj-alpha-'}
 \dim_H\F_{k,j}^-(t,\beta)\ge \dim_H\pi_{\ga_j}(\Phi_t(\F_{k,j}^-(t,\beta))).
 \end{equation}
 So it is necessary to {consider a lower bound of} $\dim_H\pi_{\ga_j}(\Phi_t(\F_{k,j}^-(t,\beta)))$.

  To do this, for each $k,j\in\N$ we {construct} a measure $\hat\mu_{k,j}$ on the tree
\[
{\T_{k,j}(t):=\set{x_1\ldots x_{n_j}^-d_1d_2\ldots: d_{nk+1}\ldots d_{nk+k}\ne (-1)^k~\forall n\ge 0}}.
\]
Let $(p_{-1}, p_0, p_1)=(\theta_k, 1-2\theta_k, \theta_k)$ be  a probability vector with $\theta_k\in[0,1/2]$ satisfying
\begin{equation}\label{eq:p-k}
\frac{1-2\theta_k}{1-(\theta_k)^k}=\beta.
\end{equation}
Then
   $\theta_k\to(1-\beta)/2$ as $k\to\f$.
 For each cylinder set of the tree $\T_{k,j}(t)$ we set
$
\hat\mu_0([x_1\ldots x_{n_j}^-])=1,
$
and {for $n\geq1$, we {let}
\[
\hat\mu_n([x_1\ldots x_{n_j}^- d_1\ldots d_{nk}])=\prod_{i=0}^{n-1}\frac{p_{d_{ik+1}}p_{d_{ik+2}}\cdots p_{d_{ik+k}}}{1-p_{-1}^k}.
\]
By} Kolmogorov's extension theorem (cf.~\cite{Durrett-1996}) there exists a unique probability measure {$\hat\mu=\hat\mu_{k,j}$}   on the tree $\T_{k,j}(t)$ satisfying
$
\hat\mu([x_1\ldots x_{n_j}^- d_1\ldots d_{nk}])=\hat\mu_n([x_1\ldots x_{n_j}^- d_1\ldots d_{nk}])$ for any $n\ge 0.
$

\begin{lemma}
  \label{lem:full-measure}
  $\hat\mu(\Phi_t(\F_{k,j}^-(t,\beta)))=1.$
\end{lemma}
\begin{proof}
  Note that $\Phi_t(\F_{k,j}^-(t,\beta))\subset \T_{k,j}(t)$. By (\ref{eq:freqency-F-kj}) it suffices to prove that for $\hat\mu$-a.e.~$(c_i)\in \T_{k,j}(t)$ we have {$freq_0((c_i))=\beta$}. Note    that $c_1\ldots c_{n_j}=x_1\ldots x_{n_j}^-$ is fixed, and the choice of each block $c_{n_j+nk+1}\ldots c_{n_j+nk+k}$ is independent and identical distributed for all $n\ge 0$ according to our definition of $\hat\mu$. So, by the law of large numbers  it follows that for $\hat\mu$-a.e.~$(c_i)\in\T_{k,j}(t)$
  \begin{align*}
  {freq_0((c_i))}&=\lim_{n\to\f}\frac{\#\set{i\in[1,n]: c_i=0}}{n}\\
  &=\lim_{n\to\f}\frac{\#\set{i\in[n_j+1, n_j+nk]: c_i=0}}{nk}\\
  &=\sum_{d_1\ldots d_k\ne (-1)^k}\frac{\#\set{i\in[1,k]: d_i=0}}{k}\hat\mu([x_1\ldots x_{n_j}d_1\ldots d_{k}])\\
  &= \sum_{\ell=1}^k\frac{\ell}{k} \binom{k}{\ell}2^{k-\ell} \frac{(\theta_k)^{k-\ell}( 1-2\theta_k )^\ell}{1-(\theta_k)^k},
  \end{align*}
  where   the last equality  follows since the number of blocks $d_1\ldots d_k\in\set{-1,0,1}^k\setminus\set{(-1)^k}$ {with} precisely $\ell(\ge 1)$ zeros is $\binom{k}{\ell}2^{k-\ell}$.
  Rearranging the above summation and by (\ref{eq:p-k}) we conclude  that
  \[
  {freq_0((c_i))}=\frac{1}{1-(\theta_k)^k}\sum_{\ell=1}^k\frac{(k-1)!}{(\ell-1)!(k-\ell)!}(2\theta_k)^{k-\ell}(1-2\theta_k)^\ell=\frac{1-2\theta_k}{1-( \theta_k )^k}=\beta,
  \]
 completing the proof.
\end{proof}
By the same argument as in the proof of Lemma \ref{lem:full-measure} one can verify that for $\hat\mu$-a.e.~$(c_i)\in\T_{k,j}(t)$ the frequencies of digits $1$ and $-1$ in $(c_i)$ are {given respectively} by
\begin{equation}
  \label{eq:freqency-mu}
  freq_{1}((c_i))=\frac{\theta_k}{1-(\theta_k)^k}\to\frac{1-\beta}{2}\quad\textrm{and}\quad freq_{-1}((c_i))=\frac{\theta_k-(\theta_k)^k}{1-(\theta_k)^k}\to\frac{1-\beta}{2}
\end{equation}
as $k\to\f$, where the limits follow {by} $\lim_{k\to\f}\theta_k=(1-\beta)/2$.

\begin{proof}
  [Proof of Proposition \ref{th:level-set} ({a lower bound})]
Note by Lemma \ref{lem:full-measure} the set $\pi_{\ga_j}(\Phi_t(\F_{k,j}^-(t,\beta)))$ has full     $\mu:= \hat\mu\circ\pi_{\ga_j}^{-1}$ measure. Then for $\mu$-a.e.~$y=\pi_{\ga_j}(x_1\ldots x_{n_j}^-d_1d_2\ldots)\in\pi_{\ga_j}(\Phi_t(\F_{k,j}^-(t,\beta)))$
\begin{align*}
 \liminf_{r\to 0^+}\frac{\log\mu(B(y,r))}{\log r}&=\liminf_{n\to\f}\frac{\log\hat\mu([x_1\ldots x_{n_j}^- d_1\ldots d_{nk}])}{\log \ga_j^{n_j+nk}}\\
 &=\liminf_{n\to\f}\frac{\log\prod_{i=1}^{nk}p_{d_i}-\log(1-p_{-1}^k)^n}{(n_j+nk)\log\ga_j}\\
 &=\liminf_{n\to\f}\frac{\sum_{s=-1,0,1}N_{s}(nk)\log {p_s}-n\log(1-p_{-1}^k)}{nk\log\ga_j},
\end{align*}
where $N_s(nk)$ denotes the number of digit $s$ in the block $d_1\ldots d_{nk}$. Therefore, by Lemma \ref{lem:full-measure} and (\ref{eq:freqency-mu}) it follows that for $\mu$-a.e.~$y \in\pi_{\ga_j}(\Phi_t(\F_{k,j}^-(t,\beta)))$
\begin{align*}
 \liminf_{r\to 0^+}\frac{\log\mu(B(y,r))}{\log r}&=\frac{\frac{\theta_k-(\theta_k)^k}{1-(\theta_k)^k}\log\theta_k+\beta\log(1-2\theta_k)+\frac{\theta_k}{1-(\theta_k)^k}\log\theta_k}{\log\ga_j}
 {-\frac{\log(1-(\theta_k)^k)}{k\log\ga_j}}\\
 &\to\quad \frac{\frac{1-\beta}{2}\log\frac{1-\beta}{2}+\beta\log\beta+\frac{1-\beta}{2}\log\frac{1-\beta}{2}}{\log\ga_j}=\frac{h(\frac{1-\beta}{2}, \beta, \frac{1-\beta}{2})}{-\log\ga_j}
\end{align*} as $k\to\f$. {Then} by Billingsley's Lemma (cf.~\cite{Falconer_1997}) it follows that
\[
\dim_H\pi_{\ga_j}(\Phi_t(\F_{k,j}^-(t,\beta)))\ge \frac{h(\frac{1-\beta}{2}, \beta, \frac{1-\beta}{2})}{-\log\ga_j}\quad {\forall j\geq1.}
\]
{This, together with (\ref{eq:subset-F-kj-alpha-}) and (\ref{eq:subset-F-kj-alpha-'}), implies}
\[
\dim_H(\La_\beta(t)\cap(\ga_j, 1/3))\ge \frac{h(\frac{1-\beta}{2}, \beta, \frac{1-\beta}{2})}{-\log\ga_j}\quad\forall j\ge 1.
\]
Since $\ga_j\nearrow 1/3$ as $j\to\f$, {we conclude that }
$
\dim_H\La_\beta(t)\ge {h(\frac{1-\beta}{2}, \beta, \frac{1-\beta}{2})}/{\log 3}.
$
\end{proof}

{Next we consider the upper bound.}
\begin{lemma}
  \label{lem:upper-alpha-4/27-1}
  Let $t\in[4/27, 1/3)\cup(1/3,1)$ and $\beta\in[0,1]$. Then
  \[
  \dim_P\La_\beta(t)\le\frac{h(\frac{1-\beta}{2}, \beta, \frac{1-\beta}{2})}{\log 3}.
  \]
\end{lemma}
\begin{proof}
{By the same argument as in the proof of Lemma \ref{lem:upper-local-lambda-4/27-1} it follows that for any $\theta\in(t,1/3)$ there exists $C>0$, such that $|\pi_\theta(\Phi_t(\lambda_1))-\pi_\theta(\Phi_t(\lambda_2))|\geq C|\lambda_1-\lambda_2|$ for any $\lambda_1,\lambda_2\in\Lambda_\beta(t)\cap(0,\theta]$. This implies that
\[\dim_P(\Lambda_\beta(t)\cap(0,\theta])\leq \dim_P {\pi_\theta} (\Phi_\theta(\Lambda_\beta (t)))\leq \dim_P\pi_{1/3}(\Phi_t(\Lambda_\beta(t))).\]
By the countable stability of packing dimension we obtain }
\begin{equation}\label{eq:upper-alpha-2}
\dim_P\La_\beta(t) \le\dim_P\pi_{1/3}(\Phi_t(\La_\beta(t))).
\end{equation}

Define a {Bernoulli} measure $\hat\nu$ on the {symbolic space $\{-1,0,1\}^{\mathbb{N}}$} such that
\[
\hat\nu([d_1\ldots d_n])=\prod_{i=1}^n p_{d_i}\quad\forall n\ge 0,
\]
where $p_{-1}=p_1=(1-\beta)/2$ and $p_0=\beta$. Let $\nu=\hat\nu\circ\pi_{1/3}^{-1}$. Note that {for each $(d_i)\in\Phi_t(\La_\beta(t))$ we have $freq_0((d_i))=\beta$}.  Then for any $y\in\pi_{1/3}(\Phi_t(\Lambda_\beta(t)))$ with {$(c_i)=\Phi_t(1/3)$}
\begin{align*}
  \limsup_{r\to 0^+}\frac{\log\nu(B(y,r))}{\log r}&=\limsup_{n\to\f}\frac{\log\prod_{i=1}^n p_{c_i}}{-\log3^{N+n}}=\frac{1}{-\log 3}\limsup_{n\to\f}\frac{\sum_{i=1}^n \log p_{c_i}}{n}\\
  &=\frac{1}{-\log 3}\limsup_{n\to\f}\frac{1}{n}\left(N_0(n)\log\beta+(n-N_0(n))\log\frac{1-\beta}{2}\right)\\
  &=\frac{1}{-\log 3}\left(\beta\log\beta+(1-\beta)\log\frac{1-\beta}{2}\right)=\frac{h(\frac{1-\beta}{2}, \beta, \frac{1-\beta}{2})}{\log 3},
\end{align*}
where $N_0(n)$ denotes the number of digit zero in the block $c_1\ldots c_n$. So, by {\cite[Proposition 2.3]{Falconer_1997}}   and (\ref{eq:upper-alpha-2}) it follows that
\[
\dim_P\La_\beta(t) \le\dim_P\pi_{1/3}(\Phi_t(\La_\beta(t)))\le \frac{h(\frac{1-\beta}{2}, \beta, \frac{1-\beta}{2})}{\log 3},
\]
completing the proof.
\end{proof}

\begin{proof}
  [Proof of Proposition \ref{th:level-set} ({an upper bound})]
  By Lemma \ref{lem:upper-alpha-4/27-1} it suffices to consider $t\in(0,4/27)$. If $t\in(0,1/9]$, then by Lemma  \ref{lem:upper-local-lambda-0-1/9} and the same argument as in the proof of Lemma \ref{lem:upper-alpha-4/27-1} it follows that
  \begin{equation*}
    \label{eq:upper-alpha-3}
    \begin{split}
      \dim_P(\La_\beta(t)\cap(0,\la_\diamond])&\le \frac{h(\frac{1-\beta}{2}, \beta,\frac{1-\beta}{2})}{-\log\la_\diamond}\le \frac{h(\frac{1-\beta}{2}, \beta,\frac{1-\beta}{2})}{\log 3},\\
      \dim_P(\La_\beta(t)\cap(\la_\diamond,\sqrt{t}))&\le \frac{h(\frac{1-\beta}{2}, \beta,\frac{1-\beta}{2})}{-\log\sqrt{t} }\le \frac{h(\frac{1-\beta}{2}, \beta,\frac{1-\beta}{2})}{\log 3}\\
      \dim_P(\La_\beta(t)\cap[\sqrt{t}, 1/3])&\le\frac{h(\frac{1-\beta}{2}, \beta,\frac{1-\beta}{2})}{\log 3}.
    \end{split}
  \end{equation*}
This implies
  $
  \dim_P\La_\beta(t)\le  {h(\frac{1-\beta}{2}, \beta,\frac{1-\beta}{2})}/{\log 3}
  $ as desired.

  Similarly, if $t\in(1/9, 4/27)$, then by Lemma \ref{lem:upper-local-lambda-1/9-4/27} and the same argument as in the proof of Lemma \ref{lem:upper-alpha-4/27-1} we can also prove $
  \dim_P\La_\beta(t)\le h(\frac{1-\beta}{2}, \beta,\frac{1-\beta}{2})/\log 3.
 $
  \end{proof}

\begin{proof}
  [Proof of  Theorem \ref{th:La(t)-alpha}] By Propositions \ref{prop:Lambda-not} and \ref{th:level-set} it suffices to prove that for any $\beta\in[0,1]$ both $\La_{not}(t)$ and $\La_\beta(t)$ are  dense in $\La(t)$. Since the proofs are similar, we only prove it for $\La_\beta(t)$. Take $\la_*\in\La(t)$ and $\de>0$.  Then we can always find $\la\in\La(t)$ such that $\Phi_t(\la)$ has a long common prefix with $\Phi_t(\la_*)$, and the tail sequence of $\Phi_t(\la)$ has digit zero frequency equaling $\beta$. In other words, $\la\in\La_\beta(t)\cap(\la_*-\de, \la_*+\de)$. So,    {$\La_\beta(t)$ is dense in  $\La(t)$.}
\end{proof}

%
%

%
%

{\section*{Acknowledgements}
The second author was supported by NSFC No.~11971079.


\end{document}